\documentclass[aop,MSNbibl,dvips]{arximspdf}
\usepackage{url,breakurl}

\doi{10.1214/12-AOP822} 
\volume{42}
\issue{2}
\pubyear{2014}
\firstpage{464}
\lastpage{496}

\makeatletter
\newcommand{\rrvert}{\vert}
\newcommand{\llvert}{\vert}

\newtheorem{theorem}{Theorem}[section]
\newtheorem{prop}[theorem]{Proposition}
\newtheorem{question}{Question}[section]
\newtheorem{lemma}[theorem]{Lemma}
\newtheorem{cor}[theorem]{Corollary}
\newtheorem{claim}[theorem]{Claim}

\newproclaim{remark}{Remark}

\def\P{\mathbb{P}}
\newcommand{\one}{\mathbf{1}}
\newcommand{\deq}{\stackrel{\triangle}{=}}
\newcommand{\reff}{R_{\mathrm{eff}}}

\renewcommand{\epsilon}{\varepsilon}
\newcommand{\var}{\operatorname{Var}}

\newcommand{\N}{{\mathbb N}}
\newcommand{\E}{{\mathbb E}}

\makeatother

\begin{document}
\begin{frontmatter}

\title{Asymptotics of cover times via Gaussian free fields:
Bounded-degree graphs and general trees}
\runtitle{Asymptotics of cover times via Gaussian free fields}

\begin{aug}
\author[A]{\fnms{Jian} \snm{Ding}\corref{}\ead[label=e1]{jianding@galton.uchicago.edu}}
\runauthor{J. Ding}
\affiliation{Stanford University and University of Chicago}
\address[A]{Department of Statistics\\
University of Chicago\\
Chicago, Illinois 60637\\
USA\\
\printead{e1}} 
\end{aug}

\received{\smonth{12} \syear{2011}}
\revised{\smonth{11} \syear{2012}}

%
\begin{abstract}
In this paper we show that on bounded degree graphs and general trees,
the cover time of the simple random walk is asymptotically equal to the
product of the number of edges and the square of the expected supremum
of the Gaussian free field on the graph, assuming that the maximal
hitting time is significantly smaller than the cover time. Previously,
this was only proved for regular trees and the 2D lattice. Furthermore,
for general trees, we derive exponential concentration for the cover
time, which implies that the standard deviation of the cover time is
bounded by the geometric mean of the cover time and the maximal hitting
time.
\end{abstract}

%

\end{frontmatter}

\section{Introduction}

Consider a random walk on a finite connected graph $G = (V,
E)$, and let $\tau_{\mathrm{cov}}(G)$ be the stopping time when the
random walk has visited every vertex in the graph for the first
time. The following fundamental parameter is known as the
\textit{cover time}:
\[
t_{\mathrm{cov}}(G) = \max_{v\in V}\E_v
\tau_{\mathrm{cov}}(G).
\]
In addition, let
$t_{\mathrm{hit}}(u, v; G)$ be the expected time it takes the random
walk started at the vertex $u$ to hit the vertex $v$, and define the
maximal hitting time $t_{\mathrm{hit}}(G) = \max_{u, v}
t_{\mathrm{hit}}(u, v; G)$. In this paper we investigate the
asymptotic value of the cover time for bounded degree graphs and
general trees as $|V| \to\infty$, and strengthen a connection between
the cover time and
the Gaussian free field.

Recall that a Gaussian free field (GFF) on the graph $G$ is a
centered Gaussian process $\{\eta_v\}_{v\in V}$ with $\eta_{v_0} =
0$ for some fixed $v_0 \in V $, and the process is characterized by the
relation $\E(\eta_u - \eta_v)^2 = R_{\mathrm{eff}}(u, v)$ for all
$u, v \in V$, where $R_{\mathrm{eff}}$ denotes the effective
resistance on $G$; see Section \ref{secprelim}. We are now ready to
state our main
results:

\begin{theorem}\label{thm-asympt}
Consider a sequence of graphs $G_n = (V_n, E_n)$ with maximal degree
bounded by a fixed $\Delta>0$ such that $t_{\mathrm{hit}}(G_n) =
o(t_{\mathrm{cov}}(G_n))$ as $n\to\infty$. For each $n$, let
$\{\eta_v\}_{v\in V_n}$ be a Gaussian free field on $G_n$ with
$\eta_{v_0^n} = 0$ for a certain $v_0^n \in V_n$. Then as $n\to
\infty$, we have
%
\begin{equation}
\label{eq-cover-asymp} t_{\mathrm{cov}}(G_n) = \bigl(1+o(1)
\bigr)|E_n| \Bigl(\E\sup_{v\in V_n} \eta_v
\Bigr)^2.
\end{equation}
\end{theorem}

\begin{remark*} Note that the expectation of the supremum for a
Gaussian free field does not depend on the choice of $v_0$, since
selecting a different ``$v_0$'' corresponds to merely shifting the
whole process by an additive mean-zero Gaussian variable.

The asymptotic identity (\ref{eq-cover-asymp}) arose in a recent
work of Ding, Lee and Peres \cite{DLP10}, where a useful connection
between cover times, Gaussian processes and Fernique--Talagrand
majorizing measure theory \cite{Fernique74,Talagrand87,Talagrand05} was
discovered. It was shown that the cover time of
any graph $G$ is equivalent to the product of the number of edges
and the square of the expected maximum of the GFF, up to a universal
multiplicative constant. In particular, the upper bound in
(\ref{eq-cover-asymp}) was established. This led to a deterministic
polynomial-time algorithm to approximate the cover time up to a
constant, which improved upon the $O((\log\log n)^2)$-approximation
for $n$-vertex graphs due to Kahn et al.
\cite{KKLV00}, and resolved a question due to Aldous and Fill
\cite{AF}.

Theorem \ref{thm-asympt} sharpens the above-mentioned universal
constant to 1 for bounded degree graphs, under the assumption that
the maximal hitting time is of smaller order than the cover time.
Two nontrivial graphs for which (\ref{eq-cover-asymp}) has been
verified are regular trees by Aldous \cite{Aldous91} and 2D lattices
by Bolthausen, Deuschel and Giacomin \cite{BDG01} and Dembo et al.
\cite{DPRZ04}. The asymptotic identity
(\ref{eq-cover-asymp}) suggests a fundamental connection between
cover times and Gaussian free fields. In addition, it may be a
useful step toward approximating the cover time algorithmically up
to a factor of $(1+\epsilon)$, where $\epsilon>0$ is arbitrarily
small.
\end{remark*}

\begin{remark*} After the current work was posted on arXiv, a
deterministic PTAS for Computing the Supremum of Gaussian Processes
was found by Meka \cite{Meka12}. Combined with our result, this gives
a deterministic PTAS for computing the cover time on bounded-degree
graphs where the maximal hitting times are significantly smaller than
the cover times.

For cover times on trees, we obtain the following exponential
concentration.
\end{remark*}

\begin{theorem}\label{thm-concentration-cover-tree}
Consider a tree $T = (V, E)$ with root $v_0\in V$. Denote by $R$ the
diameter of $T$. Let $\{\eta_v\}_{v\in V}$ be a Gaussian free field
on $T$ with $\eta_{v_0} = 0$. Then for the random walk started at
$v_0$ and any $\lambda\geq1$,
\[
\P\Bigl(\Bigl\llvert\tau_{\mathrm{cov}}(T) - |E| \Bigl(\E\sup
_v \eta_v\Bigr)^2\Bigr\rrvert\geq
\lambda|E| \sqrt{R} \E\sup_v\eta_v \Bigr)\leq C
\mathrm{e}^{-c \lambda},
\]
where $C, c>0$ are universal constants.
\end{theorem}

The following well-known commute time identity \cite{CRRST96} gives a
useful connection between random walk and \textit{effective resistance distance}:
%
\begin{equation}
\label{eq-commute-time} \kappa(u, v) = 2 |E| R_{\mathrm{eff}}(u, v),
\end{equation}
where $\kappa(u, v)$ is the commute time between $u$ and $v$ (i.e.,
the expected time it takes the random walk to travel from $u$ to $v$
and then return to $u$), and $R_{\mathrm{eff}}(u, v)$ is the effective
resistance between $u$ and $v$; see Section \ref{secprelim} for
background on electric networks. In the particular cases for trees, the
commute time identity yields that $t_{\mathrm{hit}}(T) \geq|E| R$.
Together with the result from \cite{DLP10} that $t_{\mathrm{cov}}\geq
c |E| \E(\sup_v \eta_v)^2$ for an absolute constant $c>0$, we see
from Theorem \ref{thm-concentration-cover-tree} that
%
\begin{eqnarray}
\label{eq-added} t_{\mathrm{cov}}(T) &=& |E| \Bigl(\E\sup_v
\eta_v\Bigr)^2 + O(1) |E| \sqrt{R} \E\sup
_v \eta_v \nonumber\\[-8pt]\\[-8pt]
&=& |E| \Bigl(\E\sup
_v \eta_v\Bigr)^2 + O(1)
\sqrt{t_{\mathrm{cov}}(T) t_{\mathrm
{hit}}(T)}.\nonumber
\end{eqnarray}
Now the following corollary is obvious from Theorem
\ref{thm-concentration-cover-tree} and (\ref{eq-added}).
%
\begin{cor}\label{cor-concentration-cover-tree}
Consider a sequence of trees $T_n = (V_n, E_n)$ with root \mbox{$v_0^n\in
V_n$}. For each $n$, let $\{\eta_v\}_{v\in V_n}$ be a Gaussian free
field on $T_n$ with $\eta_{v_0^n} = 0$ for all $n\in\N$. Then there
exist universal constants $c, C>0$ such that for any $\lambda\geq
0$,
\[
\P\bigl(\bigl|\tau_{\mathrm{cov}}(T_n) - t_{\mathrm{cov}}(T_n)\bigr|
\geq\lambda\sqrt{t_{\mathrm{cov}}(T_n) \cdot t_{\mathrm{hit}}(T_n)}
\bigr)\leq C \mathrm{e}^{-c\lambda}.
\]
Assume in addition $t_{\mathrm{hit}}(T_n) =
o(t_{\mathrm{cov}}(T_n))$. Then
\[
t_{\mathrm{cov}}(T_n) = \bigl(1+o(1)\bigr) \cdot|E_n|
\Bigl(\E\sup_{v\in
V_n}\eta_v \Bigr)^2.
\]
\end{cor}

For convenience, we will work exclusively with \textit{continuous-time}
Markov chains, where the transition rates between nodes are given by
the probabilities $p_{xy}$ from the discrete chain. One way to realize
the continuous-time chain is by making jumps according to the
discrete-time chain, where the times spent between jumps are i.i.d.
exponential random variables with mean 1; see \cite{AF}, Chapter 2, for
background and relevant definitions.

Note that our results automatically extend to discrete time
random walk. Let $\tau_{\mathrm{cov}}^\star$ be the cover time for
the discrete time random walk. It is clear that
$\E_v\tau_{\mathrm{cov}}^\star= \E_v \tau_{\mathrm{cov}}$ for all
$v\in V$, and therefore Theorem \ref{thm-asympt} extends to discrete
case trivially. Furthermore, the number of steps $N(t)$ performed by
a continuous-time random walk up to time $t$, has Poisson
distribution with mean $t$. Therefore, $N(t)$ exhibits a Gaussian-type
concentration around $t$ with standard deviation bounded by
$\sqrt{t}$. This implies that the concentration result in
Theorem \ref{thm-concentration-cover-tree} holds for discrete-time
case.

We remark that the assumption $t_{\mathrm{hit}} =
o(t_{\mathrm{cov}})$ is very natural. For one thing, without this
assumption, the asymptotic identity is not necessarily true. In the
case of a line on $n$ vertices, it is clear that $t_{\mathrm{cov}} =
(1+o(1))5n^2/4$ since the worst starting point is the middle point of
the line; one could see, for example, \cite{Durrett}, Exercise~4.7.3,
for estimates on expected hitting times of 1D simple random walk, while
$\E\sup_v \eta_v = (1+o(1))\sqrt{2n/\pi}$ for the
GFF $\{\eta_v\}$ (this can be deduced by the fact that $\sup_v
\eta_v$ is asymptotical to the supremum of a Brownian motion, which
has the distribution as the absolute value of a Gaussian variable. One
could also see, e.g., \cite{Durrett}, Example 7.4.3).
For another, the asymptotics of the expectation is of most interest
when the cover time is concentrated around the expectation (i.e., $\tau
_{\mathrm{cov}}/t_{\mathrm{cov}}$ converges to 1 with probability
tending to 1 as $|V| \to\infty$). It turns out
the ratio between the maximal hitting time and cover time governs the
concentration property of the cover time.
If the maximal hitting time and cover time have the
same order, it was shown that $\tau_{\mathrm{cov}}$ is not
concentrated by Aldous \cite{Aldous91b}, Proposition~1. However, $\tau
_{\mathrm{cov}}$ does exhibit
concentration around its mean under the assumption $t_{\mathrm{hit}}
= o(t_{\mathrm{cov}})$, due to the following result proved in \cite
{Aldous91b}.

\begin{theorem}[(\cite{Aldous91b})]\label{thm-Aldous}
Consider a sequence of graphs $G_n = (V_n, E_n)$ such that
$t_{\mathrm{hit}}(G_n) = o(t_{\mathrm{cov}}(G_n))$. Then with high
probability,
\[
\tau_{\mathrm{cov}}(G_n) = \bigl(1+o(1)\bigr) t_{\mathrm{cov}}(G_n).
\]
\end{theorem}

It is interesting to study the concentration of
$\tau_{\mathrm{cov}}$ quantitatively. Our
Theorem~\ref{thm-concentration-cover-tree} follows in this line of
research when the underlying graph is a general tree. In particular,
Theorem \ref{thm-concentration-cover-tree} proves that
$\tau_{\mathrm{cov}}$ exhibits an exponential concentration (which
was observed for the supremum of a Gaussian process (see, e.g.,
\cite{Ledoux89}, Theorem 7.1, Equation (7.4)) and its standard
deviation is bounded from
above by the geometric mean of the maximal hitting time and the
cover time. This seems to be the first exponential concentration
result of this type to our knowledge.

As mentioned earlier, the upper bound for (\ref{eq-cover-asymp}) has
been established in~\cite{DLP10}.

\begin{prop}[(\cite{DLP10})]\label{prop-DLP}
There exists a universal constant $C>0$, such that for any graph $G=
(V, E)$ with $v_0\in V$, we have
\[
t_{\mathrm{cov}} \leq\biggl(1+ C\sqrt{\frac{t_{\mathrm
{hit}}}{t_{\mathrm{cov}}}} \biggr) \cdot|E| \cdot
\Bigl(\E\sup_{v\in V} \eta_v \Bigr)^2,
\]
where $\{\eta_v\}_{v\in V}$ is the Gaussian free field on graph $G$
with $\eta_{v_0} = 0$.
\end{prop}
The lower bound for the cover time seems to be much more elusive.
Indeed, most of the work in \cite{DLP10} was devoted to prove that
the cover time is bounded from below via the GFF up to a universal
constant for any graph. Sharpening such constant is significantly
more challenging, partly because a fundamental ingredient of
\cite{DLP10}, known as the majorizing measure theory, loses a
multiplicative constant to begin with. As a preliminary (but
important) step to approach the lower bound, we relax the problem
based on Theorem \ref{thm-Aldous}.

\begin{theorem}\label{thm-tau-cov}
Consider a graph $G = (V, E)$ with maximal degree bounded by a fixed
$\Delta>0$. Let $\{\eta_v\}_{v\in V}$ be a Gaussian free field on
$G$ with $\eta_{v_0} = 0$ for a certain $v_0\in V$. Fix any
$0<\epsilon\leq1/10$, and assume that
%
\begin{equation}
\label{eq-assumption}t_{\mathrm{hit}} \leq\frac{\epsilon^4}{10^4\Delta
^2 (C \vee1)^2}t_{\mathrm{cov}},
\end{equation}
where $C$ is the universal constant in Proposition \ref{prop-DLP}.
Then there exists
$\delta= \delta(\epsilon, \Delta)
> 0$ such that
%
\begin{equation}
\label{eq-thm-6} \P\Bigl(\tau_{\mathrm{cov}} \geq(1-\epsilon) |E| \Bigl
(\E\sup
_v \eta_v\Bigr)^2\Bigr) \geq\delta.
\end{equation}
\end{theorem}

We now deduce the lower bound for (\ref{eq-cover-asymp}) from
Theorems \ref{thm-Aldous} and \ref{thm-tau-cov}. By Theorem
\ref{thm-Aldous}, we have that
for any $\epsilon>0$
\[
\P\bigl(\tau_{\mathrm{cov}}(G_n) \leq(1 + \epsilon) t_{\mathrm{cov}}
(G_n)\bigr) \to1\qquad \mbox{as } n \to\infty.
\]
Combined with (\ref{eq-thm-6}), it follows that for any $\epsilon>
0$ and sufficiently large $n$,
\[
(1 + \epsilon) t_{\mathrm{cov}} (G_n) \geq(1-\epsilon)
|E_n| \Bigl(\E\sup_v \eta_v
\Bigr)^2,
\]
which proves
the lower bound for (\ref{eq-cover-asymp}) by sending $\epsilon\to0$.

Next, we describe the main strategy to prove Theorem
\ref{thm-tau-cov}. Our proof employs the \textit{sprinkling} method as the
roadmap. This type of perturbation method was used by Ajtai,
Koml\'os and Szemer\'edi \cite{AKS82} in the study of percolation,
and found its applications later in that area; see, for example,
\cite{ABS04,BNP09}. In the setting of cover time, our main
intuition is the following: if there exists a \textit{thin point} at
time $\tau(t)$ [i.e., a~vertex which was visited by the random walk only
for a few number of times up to time $\tau(t)$], there should be a
positive chance that the random walk did not yet visit the thin
point up to time $\tau((1-\epsilon)t)$ (the sprinkling), and
therefore did
not yet cover the graph. Most of the work is then devoted to show the
existence of a thin
point.

\subsection{Related work}

There is a long history of the study of cover times, in probability,
combinatorics and theoretical computer science. We will review only
the work related to the very precise\vadjust{\goodbreak} estimates for cover times (and
the related supremum of Gaussian free field). We refer to the books
\cite{AF,LPW09} and the survey \cite{Lov96} for relevant background
material. For a more up-to-date account on the history for cover
times, see the introduction in \cite{DLP10} as well as the
references therein.

Previous to our work, the only nontrivial examples for which
(\ref{eq-cover-asymp}) has been verified are regular trees and the
2D torus. For regular trees, the asymptotics of cover times was
shown in \cite{Aldous91}, while the supremum of the Gaussian free
field was known as a folklore, and a precise estimate up to an
additive constant can be deduced by adapting Bramson's methods on
the maximal displacement of branching Brownian motion
\cite{Bramson78}. Indeed, an analogue of Bramson's result for a wide
range of branching random walks was proved by Addario-Berry and Reed
\cite{BR09}. For the 2D lattice, the asymptotics of the supremum of
the Gaussian free field was determined in \cite{BDG01}, and the
asymptotics of cover times was established in \cite{DPRZ04}. We
emphasize that in both cases, the asymptotics of cover times was
very tricky, despite the fact that the supremum of the GFFs had been
established.

There are additional high-precision estimates for cover times and
Gaussian free fields on trees and 2D lattice: for regular binary
trees $T_n$ of height $n$, Bramson and Zeitouni \cite{BZ09} proved
that $\sqrt{\tau_{\mathrm{cov}}(T_n)/2^n}$ is tight after proper
centering. For general trees, Feige and Zeitouni \cite{FZ09} studied
the computational perspective and designed a deterministic
polynomial-time algorithm to approximate the cover time up to a
factor of $(1+\epsilon)$ for any fixed $\epsilon>0$. For the 2D
lattice, in a recent breakthrough paper of Bramson and Zeitouni
\cite{BZ10}, it was shown that the supremum of the Gaussian free
field is tight after proper centering, together with an estimate on
its expectation up to an additive constant. It improved upon the
tightness result along a subsequence by Bolthausen, Deuschel and
Zeitouni \cite{BDZ10} and a super-concentration result due to
Chatterjee \cite{Chatterjee08}.

Miller and Peres \cite{MP09} studied the connection between cover
times and the mixing times for the random walks on corresponding
lamplighter graphs. In particular, they designed a procedure which
allowed them to compute the cover time up to $1+o(1)$ for a family of
graphs that satisfy some ``transient'' condition. Miller pointed out
that this procedure should also allow one to compute the supremum of
Gaussian free field up to $1+o(1)$. However, it seems that their
method could not be extended to the case for general trees---at
least not without further substantial ingredient.

Benjamini, Gurel-Gurevich and Morris showed that for bounded degree
graphs it is exponentially unlikely to cover the graph in linear
time \cite{BGM10}. This is a different type of large-deviation
result on the cover time from the one that we prove.

In a work of Ding and Zeitouni \cite{DZ11}, the second
order term for the cover time on a binary tree was pinned down, and a
discrepancy
from the supremum of GFF was demonstrated in this scale.


\subsection{Preliminaries}\label{secprelim}

\textit{Electric networks.} A \textit{network} is a finite,
undirected graph $G=(V,E)$ (possibly with self-loops), together with a
set of nonnegative
conductances $\{c_{xy}\dvtx  x,y \in V\}$ supported exactly on the edges
of $G$, that is, $c_{xy} > 0 \iff xy \in E$. The conductances are
symmetric so that $c_{xy} = c_{yx}$ for all $x,y \in V$. We will
write $c_x = \sum_{y \in V} c_{xy}$ for the total conductance at
vertex $x$. We will often use the notation $G(V)$ for a network on
the vertex set $V$. In this case, the associated conductances are
implicit. In the few cases when there are multiple networks under
consideration simultaneously, we will use the notation
$\tilde{c}_{xy}$ to refer to the conductances in $\tilde{G}$,
correspondingly. Note that a graph $G = (V, E)$ can be viewed as a
network $G=G(V)$ even without specifying the conductances. In that
case, each edge in the graph is assigned a unit conductance except that
each self-loop is assigned conductance 2, by
convention. In particular, $c_v = d_v$, where $d_v$ is the degree of
vertex $v$.

For such a network, we can consider the canonical \textit{discrete time
random walk on $G$}, whose transition probabilities are given by
$p_{xy} = c_{xy}/c_x$ for all $x,y \in V$. It is easy to see that this
defines the transition matrix of a reversible Markov chain on~$V$, and
that every finite-state reversible Markov chain arises in this way;
see~\cite{AF}, Section 3.2. The stationary measure of a vertex is
$\pi(x) = c_x/\sum_y c_y$.

Associated to such an electrical network is the classical
quantity $\reff\dvtx  V \times V \to[0, \infty]$ which
is referred to as the \textit{effective resistance} between pairs of nodes.
Furthermore,
the effective resistances form a metric (see, e.g., \cite{KR93}) which
we call resistance
metric. We refer to \cite{LPW09}, Chapter 9, and \cite{LP}, Chapter 2,
for a discussion about the
connection between electrical networks and the corresponding random
walk. In particular, a formal definition of effective resistance can be
given using such a connection as
\[
R_{\mathrm{eff}} (u, v) = \frac{1}{c_u \P(u\to v)},
\]
where $\P(u\to v)$ is the probability for a random walk started at $u$
to hit $v$ before returning to $u$.
In fact the effective resistance can be extended to $\reff\dvtx  V \times
2^V \to[0, \infty]$ (and indeed even as a function on $2^V \times
2^V$) such that
\[
R_{\mathrm{eff}}(u, S) = \frac{1}{c_u \P(u\to S)},
\]
where $\P(u\to S)$ is the probability for a random walk started at $u$
to hit $S$ before returning to $u$.


\textit{Gaussian free field.}
Consider a connected network $G(V)$. Fix a vertex $v_0 \in V$, and
consider the random process $\mathcal X = \{\eta_u\}_{u \in V}$,
where $\eta_{v_0}=0$, and $\mathcal X$ has density proportional to
%
\begin{equation}
\label{eqdensity} \exp\biggl(-\frac14 \sum_{u,v}
c_{uv} |\eta_u-\eta_v|^2 \biggr).
\end{equation}
The process $\mathcal X$ is called the discrete Gaussian free field
(GFF) associated with $G$. The following well-known identity relates
the GFF to the electric network (see, e.g., \cite{Janson97}, Theorem
9.20):
%
\begin{equation}
\label{eqresistGFF} \E(\eta_u-\eta_v)^2 =
\reff(u,v).
\end{equation}

\textit{Cover times and local times.}
For a connected network $G(V)$, let $(X_t)$ be a continuous-time random
walk on $G$ started at a certain $v_0\in V$. For a vertex $v \in V$
and time $t$, we define the \textit{local time $L_t^v$} by
%
\begin{equation}
\label{eqlocaltimedef} L_t^v = \frac{1}{c_v} \int
_0^t \one_{\{X_s = v\}} \,ds.
\end{equation}
It is obvious that local times are crucial in the study of cover
times, since
\[
\tau_{\mathrm{cov}} = \inf\bigl\{t>0\dvtx  L^v_t >0 \mbox{
for all } v\in V\bigr\}.
\]
To this end, it turns out that it is
convenient to decompose the random walk into excursions at $v_0 \in V$.
This motivates the following definition of the
inverse local time $\tau(t)$:
%
\begin{equation}
\label{eqinverselt} \tau(t) = \inf\bigl\{s\dvtx  L_s^{v_0} > t
\bigr\}.
\end{equation}
We study the cover time via analyzing the local time process
$\{L^v_{\tau(t)}\dvtx  v\in V\}$. In this way, we measure the cover time
in terms of $\tau(t)$ and note that the random walk is always at $v_0$
at $\tau(t)$.


\textit{Dynkin isomorphism theory.} The distribution of the
local times for a Borel right process can be fully characterized by
a certain associated Gaussian processes; results of this flavor go
by the name of \textit{Dynkin isomorphism theory}. Several versions
have been developed by Ray \cite{Ray63} and Knight \cite{Knight63},
Dynkin \cite{Dynkin83,Dynkin84}, Marcus and Rosen \cite{MR92,MR01},
Eisenbaum \cite{Eisenbaum95} and Eisenbaum et al. \cite
{EKMRS00}. In what follows, we present the second
Ray--Knight theorem in the special case of a continuous-time random
walk. It first appeared in \cite{EKMRS00}; see also Theorem 8.2.2 of
the book by Marcus and Rosen \cite{MR06} (which contains a wealth of
information on the connection between local times and Gaussian
processes). It is easy to verify that the continuous-time random
walk on a connected graph is indeed a recurrent strongly symmetric
Borel right process; see, for example, \cite{MR06} for relevant
definitions. Furthermore, in the case of random walk, the
associated Gaussian process turns out to be the GFF on the
underlying network.

\begin{theorem}[(Generalized second Ray--Knight isomorphism
theorem~\cite{EKMRS00})]\label{thmrayknight}
Consider a continuous-time random walk on graph $G = (V, E)$ from
$v_0\in V$. Let $\tau(t)$ be defined as in (\ref{eqinverselt}).
Denote by $\eta= \{\eta_x\dvtx  x\in V\}$ the GFF on $G$ with
$\eta_{v_0} = 0$. Let $\P_{v_0}$ and $\P^\eta$ be the laws of the
random walk and the GFF, respectively. Then for any $t > 0$ under the measure
$\P_{v_0} \otimes\P^\eta$,
%
\begin{equation}
\label{eqlaw} \bigl\{L_{\tau(t)}^x + \tfrac{1}{2}
\eta_x^2\dvtx  x\in V \bigr\} \stackrel{\mathrm{law}} {=} \bigl\{
\tfrac{1}{2} (\eta_x + \sqrt{2t})^2\dvtx  x\in V \bigr\}.
\end{equation}
\end{theorem}

\subsection{Outline of the paper}
Section \ref{sectrees} is devoted to the study of cover times on
general trees, and contains a coupling between local times and GFFs
on trees as a key ingredient. The coupling relies on the recursive
structure of the tree. In Section~\ref{secdetection}, we prove a
detection property
for GFF on bounded-degree graphs, based on the structure of a
sequential decomposition for Gaussian free field. This property then
translates to that of local times by Theorem \ref{thmrayknight}. In
Section \ref{secreconstruction}, we first set up a framework for
the reconstruction of random walk paths from local times and present
a connection between random walks and Eulerian circuits. Using such
a connection, we demonstrate the existence of the thin point. We
conclude the paper by discussions on future directions in
Section \ref{secdiscussion}.

\section{Concentration for cover times on general
trees}\label{sectrees} In this section, we establish a sharp
asymptotics for the cover times on trees together with an
exponential concentration around its mean, as incorporated in
Theorem \ref{thm-concentration-cover-tree}. The key of the proof is
a coupling between local times and GFFs on trees, as explored in
Section~\ref{sectree-domination}.

\subsection{Concentration for inverse local time}
Throughout the paper, we measure the cover time
$\tau_{\mathrm{cov}}$ by the inverse local time $\tau(t)$. This is
legitimate only if $\tau(t)$ is highly concentrated, which we show
in this subsection.

\begin{lemma}\label{lem-concentration-tau-t}
For a graph $G = (V, E)$ with $v_0\in V$, denote by $R$ the diameter
of the graph in the resistance metric. Let $\tau(t)$ be defined as
in (\ref{eqinverselt}), for $t>0$. Then, for any $\lambda\geq1$,
\[
\P\bigl(\bigl|\tau(t) - 2t|E|\bigr| \geq(\sqrt{ \lambda t R} + \lambda R)|E|\bigr
) \leq6
\mathrm{e}^{-\lambda/16}.
\]
\end{lemma}

\begin{remark*}
In the preceding lemma, $G$ does not have to be
a simple graph. In fact, the result can be extended to a general network.

In order to prove the above lemma, we need to use the following
concentration result on the sum of squares of Gaussian variables.
\end{remark*}

\begin{claim}\label{claim-gaussian-square}
Let $(X_1,\ldots, X_n)$ be a centered Gaussian vector such that $\E
X_i^2 \leq\sigma^2$ for all $1\leq i\leq n$. Take $a_i>0$ and write
$A = \sum_{i=1}^n a_i$. Then for any $\lambda>0$,
\[
\P\Biggl(\sum_{i=1}^n a_i
X_i^2 \geq\lambda A \sigma^2\Biggr) \leq2
\mathrm{e}^{-\lambda/4}.
\]
\end{claim}
\begin{pf}
First consider integers $k_i$, and let $k = \sum_{i=1}^n k_i$. By a
generalized H\"older inequality,
\[
\E\Biggl(\prod_{i=1}^n
X_i^{2k_i}\Biggr) \leq\prod_{i=1}^n
\bigl(\E X_i^{2k}\bigr)^{k_i/k} \leq
\sigma^{2k} (2k-1)!!,
\]
where the last transition follows from the fact that $\E Z^{2k} =
(2k-1)!!$ for a standard Gaussian variable $Z$. Therefore, we have
\begin{eqnarray*}
\E\exp\biggl(\frac{1}{4A \sigma^2} \sum_{i}a_i
X_i^2 \biggr) &=& \sum_{k=0}^\infty
\frac{1}{k!}\cdot\biggl(\frac{1}{4A
\sigma^2} \biggr)^k \cdot\E
\biggl(\sum_i a_i
X_i^2\biggr)^k \\
&\leq&\sum
_{k=0}^\infty\frac{(2k-1)!!}{k! 4^k} \leq\sum
_{k=0}^\infty\frac{1}{2^k}\leq2.
\end{eqnarray*}
Now an application of Markov's inequality completes the proof.
\end{pf}
\begin{pf*}{Proof of Lemma \ref{lem-concentration-tau-t}}
Note that $\E\eta_v^2 \leq R$ for all $v\in V$, by
(\ref{eqresistGFF}) and $\eta_{v_0} = 0$. In view of Theorem \ref
{thmrayknight}, we see
that (denoting by $d_v$ the degree of vertex $v$)
\[
\tau(t) = \sum_{v\in V} d_v
L^v_{\tau(t)} \preceq2|E|t + \sqrt{2t}\sum
_v d_v\eta_v + \sum
_v d_v \eta_v^2/2.
\]
This implies that
\begin{eqnarray*}
&&
\P\bigl(\tau(t) - 2|E|t \geq(\sqrt{ \lambda t R} + \lambda R)
|E|\bigr)\\
&&\qquad\leq\P
\biggl(\sqrt{2t} \sum_v d_v
\eta_v \geq\sqrt{\lambda tR} |E|\biggr) + \P\biggl(\sum
_v d_v\eta_v^2 \geq2
\lambda|E|R\biggr)
\\
&&\qquad\leq\mathrm{e}^{-\lambda/16} + 2\mathrm{e}^{-\lambda/16} \leq3
\mathrm{e}^{-\lambda/16},
\end{eqnarray*}
where the second inequality follows from the fact that $\sum_v
d_v\eta_v$ is a Gaussian with variance bounded by $4 R |E|^2$ as
well as Claim \ref{claim-gaussian-square}. For the lower bound,
Theorem~\ref{thmrayknight} gives that
\[
\tau(t) + \sum_v d_v
\eta_v^2/2 \succeq2t|E| + \sqrt{2t} \sum
_v d_v\eta_v.
\]
Therefore, we can deduce that
\begin{eqnarray*}
&&
\P\bigl(\tau(t) - 2 t|E| \leq- (\sqrt{\lambda t R} + \lambda R)|E|\bigr) \\
&&\qquad\leq\P
\biggl(\sqrt{2t} \sum_v d_v
\eta_v \leq-\sqrt{\lambda tR} |E|\biggr) + \P\biggl(\sum
_v d_v\eta_v^2 \geq2
\lambda|E|R\biggr) \\
&&\qquad\leq3 \mathrm{e}^{-\lambda/16},
\end{eqnarray*}
where we have used the fact that $\{-\eta_v\}$ has the same law as $\{
\eta_v\}$.
\end{pf*}

\subsection{Dominating local times by Gaussian free
fields}\label{sectree-domination}

In this subsection, we establish the following coupling between the
(square-root of) local times and GFF.

\begin{theorem}\label{thm-coupling-tree}
Given a tree $T = (V, E)$ rooted at $v_0$, consider the local time
process $\{L^v_{\tau(t)}\}_{v\in V}$ and the associated Gaussian
free field $\{\eta_v\}_{v\in V}$. For any $t>0$, we have
%
\begin{equation}
\label{eq-stochastic-domination} \min_{v\in V}\sqrt
{L^v_{\tau(t)}} \preceq\frac{1}{\sqrt{2}}\max\Bigl\{
\min_{v\in V}\eta_v + \sqrt{2t}, 0 \Bigr\}.
\end{equation}
\end{theorem}
The proof of the preceding theorem combines a coupling lemma for
random variables and the recursive structure of local times on
trees.

\textit{Gaussian, Poisson and exponential}: \textit{a coupling.} The
following identity in law involves Gaussian variables, Poisson
variables, as well as exponential variables. It can be viewed as a
preliminary version of the isomorphism theorem. We give a proof for
completeness.

\begin{lemma}\label{lem-baby-iso}
Let $X$ be a standard Gaussian variable and $Y_i$ be standard
exponential variables. Let $N$ be a Poisson random variable with
mean $\ell\geq0$. Suppose that all the variables are independent.
Then
%
\begin{equation}
\label{eq-baby-iso} \sum_{i=1}^N
Y_i + \frac{1}{2}X^2 \stackrel{\mathrm{law}} {=}
\frac{1}{2} (X +\sqrt{2\ell})^2.
\end{equation}
\end{lemma}
\begin{pf}
The proof is done by calculating the Laplace transforms. Fix an
arbitrary $\lambda> 0$. We start with the right-hand side. A
straightforward calculation yields that
\begin{eqnarray*}
\E\mathrm{e}^{-\lambda(X+\sqrt{2\ell})^2/2} &=& \frac{1}{\sqrt{2\pi}}\int
_{-\infty}^\infty
\mathrm{e}^{-x^2/2} \mathrm{e}^{-\lambda(x+ \sqrt{2\ell})^2/2} \,dx \\
&=&
\mathrm{e}^{-{\lambda\ell}/{(1+\lambda)}}
\frac{1}{\sqrt{2\pi}}\int_{-\infty}^\infty
\mathrm{e}^{-({1+\lambda})(x + \lambda
\sqrt{2\ell}/(1+\lambda))^2/{2}} \,dx
\\
&=& (1+\lambda)^{-1/2} \mathrm{e}^{-{\lambda
\ell}/{(1+\lambda)}} \frac{1}{\sqrt{2\pi}}\int
_{-\infty}^\infty\mathrm{e}^{-{x^2}/{2}} \,dx \\
&=& (1+
\lambda)^{-1/2} \mathrm{e}^{-{\lambda
\ell}/({1+\lambda})}.
\end{eqnarray*}
For the special case of $\ell= 0$, we get that
$\E\mathrm{e}^{-\lambda X^2/2} = 1/\sqrt{\lambda+ 1}$. For any
$\theta>0$, note that
\[
\E\theta^N = \mathrm{e}^{-\ell} \sum
_{k=0}^\infty\frac{(\ell
\theta)^k}{k!} = \mathrm{e}^{\ell(\theta- 1)}.
\]
Combined with the fact that $\E\mathrm{e}^{-\lambda Y_i} =
1/(1+\lambda)$, it follows that
\[
\E\mathrm{e}^{-\lambda(\sum_{i=1}^N Y_i + X^2/2)} = \frac{1}{\sqrt
{\lambda+1}} \cdot\mathrm{e}^{-{\ell\lambda}/({\lambda+1}) }.
\]
Thus, we have shown that the Laplace transforms of both sides are
equal, completing the proof.
\end{pf}

Based on Lemma \ref{lem-baby-iso}, we can derive the following
stochastic domination.

\begin{lemma}\label{lem-domination}
Let $X$ be a standard Gaussian variable and $Y_i$ be i.i.d.
standard exponential variables. Let $N$ be an independent Poisson
random variable with mean $\ell\geq0$. Then
\[
\sqrt{\sum_{i=1}^N
Y_i} \preceq\frac{1}{\sqrt{2}}\max\{ X+\sqrt{2\ell}, 0\}.
\]
\end{lemma}
\begin{pf}
Note that
\[
\P\Biggl(\sum_{i=1}^N Y_i =
0\Biggr) = \P(N=0) = \mathrm{e}^{-\ell}.
\]
Denoting by $f$ the density function of standard Gaussian variable,
we can then deduce that for any $x>0$,
\[
f\bigl(-(\sqrt{2\ell} + x)\bigr) = \frac{1}{\sqrt{2\pi}} \mathrm
{e}^{-{(\sqrt{2\ell} + x)^2}/{2}}
\leq\P\Biggl(\sum_{i=1}^N
Y_i = 0\Biggr) \cdot f(x).
\]
Integrating over both sides, we have
\[
\P(X/\sqrt{2} + \sqrt{\ell} \leq-x) \leq\P\Biggl(\sum
_{i=1}^N Y_i = 0\Biggr) \cdot\P(X
\leq- \sqrt{2}x).
\]
Together with (\ref{eq-baby-iso}), we obtain
that
%
\begin{eqnarray}
\label{eq-prob-compare}
&&\P(X/\sqrt{2} + \sqrt{\ell} \geq x) \nonumber\\
&&\qquad= \P
\bigl(|X/\sqrt{2} +
\sqrt{\ell}| \geq x\bigr) - \P(X/\sqrt{2} + \sqrt{\ell} \leq-x)
\nonumber\\[-8pt]\\[-8pt]
&&\qquad\geq \P\Biggl(\Biggl(\sum_{i=1}^N
Y_i + X^2/2\Biggr)^{1/2} \geq x \Biggr) - \P
\Biggl(\sum_{i=1}^N Y_i = 0
\Biggr) \cdot\P(X \leq- \sqrt{2}x)
\nonumber
\\
&&\qquad\geq \P\Biggl(\sqrt{\sum_{i=1}^N
Y_i} \geq x \Biggr).\nonumber
\end{eqnarray}
%
The desired stochastic domination follows directly from
(\ref{eq-prob-compare}).\vadjust{\goodbreak}
\end{pf}

\textit{Recursive structure of local times on trees.} The
following recursive construction of local times makes use of the
structure of trees.
%
\begin{lemma}\label{lem-GFF-tree}
For a tree $T = (V, E)$ rooted at $v_0\in V$, consider the local
time process $\{L^v_{\tau(t)}\}_{v\in V}$. For an arbitrary $v\in
V\setminus\{v_0\}$, denote by $T_v \subseteq T$ the subtree rooted
at $v$ and by $u$ its parent. We have
\[
\bigl(L^v_{\tau(t)} \mid L^u_{\tau(t)} =
\ell, \bigl\{L^w_{\tau(t)}\bigr\}_{w\in
T\setminus T_{v}}\bigr)
\stackrel{\mathit{law}} {=} \sum_{i=1}^N
Y_i,
\]
where $N$ is an independent Poisson variable with mean $\ell$ and
$Y_i$ are i.i.d. standard exponential variables.
\end{lemma}
\begin{pf}
A random walk on tree $T$ can be decomposed into a random walk on
$T\setminus T_v$ with excursions on $T_v$. More precisely, we first
take a random walk on $T\setminus T_v$ and then insert i.i.d.
excursions $\{\mathrm{Ex}_i\}$ at Poisson rate $1$ for all the time
the random walk
spends at $u$. Thus, the total number of excursions $N$ conditioning on
$L^u_{\tau(t)} = \ell$ and $ \{L^w_{\tau(t)}\}_{w\in
T\setminus T_{v}}$ is distributed as
\[
\bigl(N \mid L^u_{\tau(t)} = \ell\mbox{ and } \bigl
\{L^w_{\tau(t)}\bigr\}_{w\in
T\setminus T_{v}}\bigr) \stackrel{\mathrm{law}} {=}
\operatorname{Poi}(\ell),
\]
where $\operatorname{Poi}(\ell)$ denotes a Poisson variable with mean $\ell$.
Furthermore each excursion $\mathrm{Ex}_i$ starts traversing from $u$
to $v$
and performs a random walk on $T_v \cup\{u\}$ until going back to
$u$. Note that every time the random walk makes a
move from $v$, the chance for it to go back to $u$ is $1/d_v$.
Therefore, the time $Y'_i$ accumulated at $v$ at each excursion
$\mathrm{Ex}_i$ is distributed as
\[
Y'_i \stackrel{\mathrm{law}} {=} \sum
_{j=1}^{N_i} Z_i,
\]
where $Z_i$ are i.i.d. exponential variables with mean $1$ and $N_i$
is an independent geometric variable with mean $d_v$. Thus, $Y'_i \sim
\operatorname{Exp}(d_v)$ is an exponential variable with mean $d_v$, and
hence $\frac{1}{d_v}Y'_i \sim\operatorname{Exp}(1)$. Altogether, it gives that
\[
\bigl(L^v_{\tau(t)} \mid L^u_{\tau(t)} =
\ell\mbox{ and } \bigl\{L^w_{\tau
(t)}\bigr\}_{w\in
T\setminus T_{v}}
\bigr) \stackrel{\mathrm{law}} {=} \frac{1}{d_v}\sum_{i=1}^{N}Y'_i,
\]
which has the same distribution as claimed in the statement of the lemma.
\end{pf}

\begin{pf*}{Proof of Theorem \ref{thm-coupling-tree}}
Note that the GFF $\{\eta_v\}_{v\in V}$ on a tree can be constructed in
the following way. Let $\{X_e\}_{e\in E}$ be i.i.d. standard Gaussian
variable. Then for $v\in V$,
\[
\eta_v = \sum_e X_e,
\]
where the summation is over the edges in the path from the root
$v_0$ to $v$. Consider $v\in V$ with parent $u$,
and let $T_v$ be the subtree of $T$ rooted at $v$. We get that
\[
(\eta_v \mid\eta_u = x_u,
\eta_w = x_w \mbox{ for } w\in T\setminus
T_v) \stackrel{\mathrm{law}} {=} X + x_u,
\]
where $X$ is a standard Gaussian variable. The local time process can
be constructed in the same fashion by recursively exploring the local
times at vertices
away from the root. More precisely, we apply
Lemma \ref{lem-GFF-tree} and get that
\[
\bigl(L^v_{\tau(t)} \mid L^u_{\tau(t)} =
\ell_u, L^w_{\tau(t)} = \ell_w \mbox{
for } w\in T\setminus T_v\bigr) \stackrel{\mathrm{law}} {=} \sum
_{i=1}^N Y_i,
\]
where $N$ is a Poisson variable with mean $\ell_u$ and $Y_i$ are
i.i.d. standard exponential variables. If $0 = \sqrt{\ell_u} \leq
\max(\frac{x_u + \sqrt{2t}}{\sqrt{2}}, 0)$, we could couple the
rest of the process in an arbitrary way. Otherwise if $0< \sqrt{\ell
_u} \leq\frac{x_u + \sqrt{2t}}{\sqrt{2}}$, we use the
decompositions of GFF and local time process and apply
Lemma \ref{lem-domination}, and obtain that
\begin{eqnarray*}
&&\bigl(L^v_{\tau(t)} \mid L^u_{\tau(t)} =
\ell_u, L^w_{\tau(t)} = \ell_w
\mbox{ for } w\in T\setminus T_v\bigr)
\\
&&\qquad\preceq\bigl(\tfrac{1}{\sqrt{2}}\max\{\eta_v + \sqrt{2t}, 0 \}
\mid\eta_u = x_u, \eta_w = x_w
\mbox{ for } w\in T\setminus T_v\bigr).
\end{eqnarray*}
It is a well-known fact that for random variables $Z_1$ and $Z_2$, we
have $Z_1 \preceq Z_2$ if and only if there exists a coupling $(Z_1,
Z_2)$ such that
$Z_1 \leq Z_2$. Then it follows that given $\{\eta_u = x_u, \eta_w =
x_w$ for $w\in T\setminus T_v\}$ and $L^u_{\tau(t)} = \ell_u,
L^w_{\tau(t)} = \ell_w$ for $w\in T\setminus T_v)$ with
$\sqrt{\ell_u} \leq\frac{x_u + \sqrt{2t}}{\sqrt{2}}$, there
exists a coupling such that
%
\begin{equation}
\label{eq-to-couple} \sqrt{L^v_{\tau(t)}}\leq
\tfrac{1}{\sqrt{2}}\max\{\eta_v + \sqrt{2t}, 0 \}.
\end{equation}
Applying (\ref{eq-to-couple}) recursively completes the proof.
\end{pf*}

In order to establish the concentration of cover times on general
trees, we need the following classical result on the concentration
of Gaussian processes; see, for example, \cite{Ledoux89},
Theorem 7.1, Equation (7.4).

\begin{lemma}\label{lem-gaussian-concentration}
Consider a Gaussian process $\{\eta_x\dvtx  x\in V\}$, and define $\sigma
= \sup_{x \in V} (\E(\eta_x^2))^{1/2}$. Then for $\alpha>0$,
\[
\P\Bigl(\Bigl\llvert\sup_{x\in V} \eta_x - \E\sup
_{x \in V} \eta_x\Bigr\rrvert> \alpha\Bigr) \leq2
\exp\bigl(-\alpha^2 / 2\sigma^2\bigr).
\]
\end{lemma}

\begin{pf*}{Proof of Theorem \ref{thm-concentration-cover-tree}}
We first consider the upper bound on $\tau_{\mathrm{cov}}$.
Let $t^+ = (\E\sup_v \eta_v +
\beta\sqrt{R})^2/2$, for $\beta>0$ to be specified. Note that
here
$R$ is the diameter of the tree, and thus also the diameter in the
resistance metric.\vadjust{\goodbreak} Therefore, we have $\E\sup_v\eta_v \geq\sqrt{R/2\pi}$.
Observe that on the event
$\{\tau_{\mathrm{cov}}
> \tau(t^+)\}$, there exists at least one vertex $v\in V$ such that
$L^v_{\tau(t^+)} = 0$.
Fix an arbitrary ordering on $V$, and let $Z$ be the first vertex
such that $L^Z_{\tau(t)} = 0$ if $\tau_{\mathrm{cov}} > \tau(t^+)$.
Since $\E\eta_v^2 \leq R$ for all $v\in V$, we have $\P(\eta_v^2
\geq\beta^2
R/4) \leq2\mathrm{e}^{-\beta^2/8}$. Since $\{\eta_v\}_{v\in V}$ and
$\{L^v_{\tau(t^+)}\}_{v\in V}$ are two independent processes, we
obtain
%
\begin{eqnarray}
\label{eq-decompose}
&&
\P\biggl( \bigl\{\tau_{\mathrm
{cov}} > \tau\bigl(t^+\bigr)
\bigr\} \Bigm\backslash\biggl\{\exists v\in V\dvtx L^v_{\tau(t^+)} +
\frac{1}{2}\eta_v^2 < \frac{\beta^2 R}{8} \biggr\}
\biggr) \nonumber\\
&&\quad\leq\P\biggl(\eta_Z^2 \geq\frac{\beta^2 R}{4}
\Bigm|\tau_{\mathrm{cov}} > \tau\bigl(t^+\bigr)\biggr) \\
&&\quad\leq2\mathrm
{e}^{-\beta^2/8}.\nonumber
\end{eqnarray}
On the other hand, we
deduce from Lemma \ref{lem-gaussian-concentration} with $\alpha=
\beta\sqrt{R}/2$ that
\[
\P\biggl(\frac{1}{2}\inf_v \bigl(\sqrt{2t^+} +
\eta_v\bigr)^2 \leq\beta^2 R/8 \biggr) \leq\P
\Bigl( \inf_v\bigl(\sqrt{2t^+} + \eta_v\bigr)
\leq\beta\sqrt{R}/2 \Bigr) \leq2\mathrm{e}^{-\beta^2/8}.
\]
Applying Theorem \ref{thmrayknight} again and combined with
(\ref{eq-decompose}), we get that
\[
\P\bigl(\tau_{\mathrm{cov}}>\tau\bigl(t^+\bigr)\bigr) \leq4\mathrm
{e}^{-\beta^2/8}.
\]
For $\lambda\geq4$, set $\beta= \sqrt{\lambda}/4$, and therefore
\begin{eqnarray*}
&&
\P\Bigl(\tau_{\mathrm{cov}} - |E|\Bigl(\E\sup_v
\eta_v\Bigr)^2 \geq\lambda|E| \Bigl(\E\sup
_v\eta_v \sqrt{R} + R\Bigr)\Bigr)
\\
&&\quad\leq\P\bigl(\tau_{\mathrm{cov}} \geq\tau\bigl(t^+\bigr)\bigr) + \P\Bigl
(\tau
\bigl(t^+\bigr) \geq|E|\Bigl(\E\sup_v\eta_v
\Bigr)^2 + \lambda|E| \Bigl(\E\sup_v
\eta_v \sqrt{R} + R\Bigr)\Bigr)
\\
&&\quad\leq 4 \mathrm{e}^{-\lambda/128} + \P\biggl( \tau\bigl(t^+\bigr) - 2
|E| t^+
\geq\sqrt{\frac{\lambda}{4} t^+ R} |E| + \frac{\lambda}{4} R
|E|\biggr)\\
&&\quad
\leq4
\mathrm{e}^{-\lambda/128} + 6 \mathrm{e}^{-\lambda/64} \leq10
\mathrm{e}^{-\lambda/128},
\end{eqnarray*}
where we have applied Lemma \ref{lem-concentration-tau-t}.

We next turn to the lower bound on $\tau_{\mathrm{cov}}$. Let $t^- =
(\E\sup_v \eta_v
-\beta\sqrt{R})^2/2$. For $\lambda\geq4$, set $\beta= \sqrt{\lambda
}/4$. We assume that $\E\sup_v \eta_v
-\beta\sqrt{R}\geq0$ (otherwise there is nothing to prove). Applying
Theorem \ref{thm-coupling-tree} together with Lemma
\ref{lem-gaussian-concentration}, we
obtain that
\[
\P\bigl(\tau_{\mathrm{cov}} \leq\tau\bigl(t^-\bigr)\bigr) \leq\P\Bigl
(\sup
_v \eta_v \leq\E\sup_v
\eta_v - \beta\sqrt{R}\Bigr) \leq2\mathrm{e^{-\beta^2/2}}.
\]
This gives that
\begin{eqnarray*}
&&
\P\Bigl(\tau_{\mathrm{cov}} - |E|\Bigl(\E\sup_v
\eta_v\Bigr)^2 \leq-\lambda|E|\Bigl(\E\sup
_v \eta_v \sqrt{ R} + R\Bigr)\Bigr)
\\
&&\quad\leq\P\bigl(\tau_{\mathrm{cov}} \leq\tau\bigl(t^-\bigr)\bigr) + \P\Bigl
(\tau
\bigl(t^-\bigr) \leq|E|\Bigl(\E\sup_v \eta_v
\Bigr)^2 -\lambda|E|\Bigl(\E\sup_v
\eta_v \sqrt{ R} + R\Bigr)\Bigr)
\\
&&\quad\leq 2\mathrm{e}^{-\lambda/32} + \P\biggl( \tau\bigl(t^-\bigr) - 2 |E| t^-
\leq- \sqrt{\frac{\lambda}{4} t^- R} - \frac{\lambda}{4} R |E|\biggr)\\
&&\quad\leq2
\mathrm{e}^{-\lambda/32} + 6 \mathrm{e}^{-\lambda/64} \leq8\mathrm
{e}^{-\lambda/64},
\end{eqnarray*}
where we have again used Lemma \ref{lem-concentration-tau-t}. This
completes the proof of Theorem~\ref{thm-concentration-cover-tree}.
\end{pf*}

\section{Detecting Gaussian free field in a tiny
window}\label{secdetection}

In this section, we establish that for a GFF on a bounded-degree
graph, there is a nonnegligible chance to detect a vertex with
value in a small window around the median of the supremum of the
Gaussian free field. Crucially, the width of the window is
interacting with the values of the GFF on the neighborhood of
detected vertex.

\subsection{A sequential decomposition of Gaussian free field}
\label{secGFF}

We will use the following network-reduction lemma, which we believe
has been known by the community for a long time. While we did not
manage to trace the original reference, we note that it has appeared
as an exercise [Exercise 2.47(d)] in the book \cite{LP}. For a proof,
see, for example, \cite{DLP10}, Lemma 2.9.
%
\begin{lemma}\label{lem-network-reduction}
For a network $G(V)$ and a subset $\tilde{V} \subset V$, there
exists a network $\tilde{G}(\tilde{V})$ such that for all $u,v\in
\tilde{V}$, we have
\[
\tilde{c}_v = c_v \quad\mbox{and}\quad R^{\tilde{G}}_{\mathrm{eff}}(u,
v) = R_{\mathrm{eff}}(u, v).
\]
We call $\tilde{G}(\tilde{V})$ the \textup{reduced network}.
Furthermore, the projection of the random walk on $G$ to $\tilde{V}$
has the same law as the random walk on $\tilde{G}$.
\end{lemma}
%
%

Based on the preceding lemma, we can now easily deduce a sequential
decomposition of GFF, which characterizes an important facet of its
special structure as well as its interplay with electric networks
and random walks. The following lemma is well known (see \cite
{Dynkin80}, Theorem
1.2.2, and \cite{Janson97}, Theorem 9.20), and we give a
proof for completeness.
%
\begin{lemma}\label{lem-DGFF}
Let $\{\eta_v\}_{v\in V}$ be a GFF on a graph $G = (V, E)$ with $\eta
_{v_0} = 0$. For
$v_0 \in S\subset V$ and $v\in V$, let $\tau$ be the hitting time to
$S$ for
a simple random walk $X_t$ on $G$, and let $a_u = \P_v(X_\tau= u)$
for $u\in S$. Then
%
\begin{eqnarray}
\label{eq-GFF-exp}
\E\bigl(\eta_v \mid\{\eta_u\}_{u\in S}\bigr) &=&
\sum_{u\in S}a_u \eta_u,
\\
\label{eq-GFF-var}
\var\bigl(\eta_v \mid\{\eta_u\}_{u\in S}\bigr)
&=& R_{\mathrm{eff}}(v, S).
\end{eqnarray}
\end{lemma}
\begin{pf}
The lemma trivially holds for $v\in S$. Therefore, we assume in what
follows $v\notin S$.
By Lemma \ref{lem-network-reduction} and the fact that the law of the GFF
is completely determined by the resistance metric [see\vadjust{\goodbreak}
(\ref{eqresistGFF})], we see that $\{\eta_w\}_{w\in S\cup\{v\}}$
has the same law as the GFF on the reduced network $\tilde{G} =
\tilde{G}(S\cup\{v\})$. Now, fix a set of real numbers
$\{g_u\}_{u\in S}$. By the definition of GFF, we have that the
conditional density of $\eta_v$ given $\{\eta_u = g_u\}_{u\in S}$
satisfies
\begin{eqnarray*}
f\bigl(\eta_v \mid\{\eta_u = g_u
\}_{u\in S}\bigr) &\propto& \exp\biggl(-\frac12 \sum
_{u\in S} \tilde{c}_{uv} |\eta_v -
g_u|^2 \biggr) \\
&\propto&\exp\biggl(-\frac{\tilde{c}_v - \tilde{c}_{v,
v}}{2}
\biggl(\eta_v - \sum_{u\in S}
\frac{\tilde{c}_{u,
v}}{\tilde{c}_v - \tilde{c}_{v, v}}g_u \biggr)^2 \biggr),
\end{eqnarray*}
where $\tilde{c}_v = \sum_{u\in S \cup\{v\}} \tilde{c}_{u, v}$.
This implies
that conditioning on $\{\eta_u = g_u\}_{u\in S}$, we have $\eta_v$
distributed as a normal variable with mean $\sum_{u\in S}
\frac{\tilde{c}_{u, v}}{\tilde{c}_v - \tilde c_{v, v}}g_u$ and variance
$1/(\tilde{c}_v - \tilde{c}_{v, v})$. Recall that the projection of
the random walk on
$G$ has the same law as the random walk on $\tilde{G}$, and we see that
$a_u = \tilde{c}_{u, v}/(\tilde{c}_v - \tilde{c}_{v, v})$. This
verifies equality
(\ref{eq-GFF-exp}). Furthermore, since the reduced network preserves
the resistance metric on the subset, we have
\[
R_{\mathrm{eff}}(v, S) = \tilde{R}_{\mathrm{eff}}(v, S) = 1/(\tilde
{c}_v - \tilde{c}_{v, v}),
\]
completing verification of (\ref{eq-GFF-var}).
\end{pf}

\subsection{Detection of Gaussian free field}

We single out the following observation on Gaussian variables that
plays a significant role in our detection argument. Roughly
speaking, the next claim captures the typical over-shoot for a
Gaussian variable conditioning on the event that its value exceeds a
certain threshold (say,~0). As an important feature, the over-shoot
can be controlled via both the standard deviation and the mean.

\begin{claim}\label{claim-gaussian}
Let $X\sim N(-\mu, \sigma^2)$ be a Gaussian variable with $\mu\geq
0$. For any $0 \leq\epsilon\leq1$, we have
\[
\P\biggl( 0 \leq X \leq\epsilon\biggl(\sigma\wedge\frac{\sigma
^2}{\mu} \biggr)
\biggr) \geq\frac{\epsilon}{5} \cdot\P(X \geq0).
\]
\end{claim}
\begin{pf}
Denote by $f$ the density of $X$ and consider $x\geq0$. It is
straightforward to check that for any $k\in\N$,
\begin{eqnarray*}
f \biggl(x+ k\epsilon\biggl(\sigma\wedge\frac{\sigma^2}{\mu} \biggr)
\biggr) &=&
\frac{1}{\sqrt{2\pi} \sigma} \exp\biggl(-\frac
{1}{2\sigma^2} \biggl(\mu+ x+ k\epsilon
\biggl(\sigma\wedge\frac{\sigma^2}{\mu} \biggr) \biggr)^2 \biggr) \\
&\leq&
\bigl(\mathrm{e}^{-{\epsilon^2 k^2}/{2}} \vee\mathrm{e}^{-k\epsilon
} \bigr) f(x).
\end{eqnarray*}
It then follows that
\begin{eqnarray*}
\P(X\geq0) &=& \sum_{k=0}^\infty\P\biggl(k
\epsilon\biggl(\sigma\wedge\frac{\sigma^2}{\mu} \biggr) \leq X \leq
(k+1) \epsilon
\biggl(\sigma\wedge\frac{\sigma^2}{\mu} \biggr) \biggr)
\\
&\leq& \P\biggl(0\leq X\leq\epsilon\biggl(\sigma\wedge\frac{\sigma
^2}{\mu}
\biggr) \biggr) \sum_{k=0}^\infty\bigl(
\mathrm{e}^{-{\epsilon^2 k^2}/{2}} \vee\mathrm{e}^{-k\epsilon}
\bigr) \\
&\leq&
\frac{5}{\epsilon} \cdot\P\biggl(0\leq X\leq\epsilon\biggl(\sigma
\wedge
\frac{\sigma^2}{\mu} \biggr) \biggr).
\end{eqnarray*}
\upqed
\end{pf}

The main proposition in this section harnesses the preceding
observation as well as the sequential decomposition of Gaussian free
field. In light of Claim \ref{claim-gaussian}, we compare the event
for detection in a tiny window to the event of exceeding the median
of the supremum of the GFF.

\begin{prop}\label{prop-gff-detection}
Given a graph $G= (V, E)$ with maximal degree bound\-ed by $\Delta$,
let $\{\eta_v\}_{v\in V}$ be the GFF on $G$ with $\eta_{v_0} = 0$
for some $v_0\in V$. For
$v\in V$, denote by $N_v$ the set of neighbors of $v$. Then for any
$0 \leq\epsilon\leq1$ and $M> 0$,
\[
\P\biggl(\exists v\in V\dvtx  M \leq\eta_{v} \leq M + \epsilon\biggl(1
\wedge\frac{\Delta}{\sum_{u\in N_v}|M-\eta_{u}|} \biggr) \biggr) \geq
\frac{2\epsilon}{ 10^{\Delta}} \P\Bigl(\sup
_v \eta_v \geq M\Bigr).
\]
\end{prop}

\begin{pf}
Write $n= |V|$. We first specify an ordering $v_0,\ldots, v_{n-1}$
(here $v_0$ is the same $v_0$ as in the statement of the proposition)
on $V$ such that for all $0< k \leq n-1$, we have $v_k \sim v_j$ for
some $j<k$. Note that such an ordering trivially exists for any
connected graph. For easiness of notation, we write $V_k = \{v_0,\ldots, v_k\}$.

Define the event
\[
A_k \deq\bigcap_{v\in V_k}\{
\eta_v < M\}.
\]
By a standard decomposition, we have
%
\begin{eqnarray}
\label{eq-decomoposition-1}
&&
\P\Bigl(\sup_v
\eta_v \geq M \Bigr)\nonumber \\
&&\qquad= \sum_{k=1}^{n-1}
\P(A_{k-1}, \eta_{v_{k}} \geq M)
\\
&&\qquad= \sum_{k=1}^{n-1}\int
_{(-\infty, M)^k} \P(\eta_{v_{k}} \geq M \mid
\eta_{v_i} = x_i \mbox{ for } 0\leq i < k )
\mu_{k-1}(dx),\nonumber
\end{eqnarray}
where $\mu_{k-1}$ is the joint density for $\{\eta_v\}_{v\in
V_{k-1}}$. Denote by
\[
\Xi= \biggl\{\exists1\leq k < n\dvtx  M \leq\eta_{v_k} \leq M + 2\epsilon
\biggl(1 \wedge\frac{\Delta}{\sum_{u\in N_{v_k}} |M-\eta_{u}|} \biggr)
\biggr\}.
\]
We have a similar decomposition,
%
\begin{eqnarray}
\label{eq-decomoposition-2}
&&\P(\Xi) \geq\sum_{k=1}^{n-1}
\int_{(-\infty, M)^k} \P\biggl(M \leq\eta_{v_{k}} \leq M + 2
\epsilon\biggl(1 \wedge\frac{\Delta}{\sum_{u\in N_{v_k}}
|M-\eta_{u}|} \biggr) \Bigm|\nonumber\\
&&\hspace*{244pt} \{\eta_{v_i} =
x_i\}_{0\leq i<k} \biggr)\\
&&\qquad\quad\hspace*{61pt}{}\times \mu_{k-1}(dx).\nonumber
\end{eqnarray}

The key to the proof of the proposition lies in a comparison for
integrands in the two decompositions. Take any $k$ and $(x_0,\ldots, x_{k-1}) \in(-\infty, M)^k$ with $x_0=0$. Write $a_{k, i} =
\P_{v_k}(\tau_{V_{k-1}} = \tau_{v_i})$ for
$i<k$, where $\tau_A$ is the hitting time to $A$ for any $A\subset
V$. By (\ref{eq-GFF-exp}),
\[
\E(\eta_{v_{k}} - M \mid\eta_{v_i} = x_i \mbox{
for } 0\leq i< k) = \sum_{i=0}^{k-1}
a_{k, i} x_i - M = -\sum_{i=0}^{k-1}
a_{k, i}| M - x_i| \leq0.
\]
By
(\ref{eq-GFF-var}), we have
\[
\var(\eta_{v_{k}} \mid\eta_{v_i} = x_i \mbox{ for
} 0\leq i < k) = R_{\mathrm{eff}}(v_{k}, V_{k-1}) \deq
\sigma_k^2 \in[1/\Delta, 1].
\]
Applying Claim \ref{claim-gaussian} and using the fact that the
conditional law of $\eta_v$ is Gaussian, we get that
%
\begin{eqnarray}
\label{eq-control-v-k}\quad
&&
\P \biggl( M \leq\eta_{v_{k}} \leq M + \epsilon
\biggl(\sigma_k \wedge\frac{\sigma_k^2}{\sum_{i<k}a_{k, i}|M-\eta
_{v_i}|} \biggr) \Bigm|
\eta_{v_i} = x_i \mbox{ for } 0\leq i < k \biggr)
\nonumber\\[-8pt]\\[-8pt]
&&\qquad\geq\frac{\epsilon}{5} \cdot\P(\eta_{v_{k}} \geq M \mid
\eta_{v_i} = x_i \mbox{ for } 0\leq i < k).\nonumber
\end{eqnarray}

We next turn to control the GFF over the neighbors of $v_k$.
Consider $x_k$ such that
%
\begin{equation}
\label{eq-x-k} M \leq x_k \leq M + \epsilon\biggl(
\sigma_k \wedge\frac{\sigma_k^2}{\sum_{i<k}a_{k, i}|M-
x_i|} \biggr).
\end{equation}
Write $f_w = f_w(x_0,\ldots, x_{k-1})
= \E(\eta_w \mid\{\eta_{v_i} = x_i\}_{0\leq i<k})$, for $w\in V$.
We claim that
%
\begin{equation}
\label{eq-control-neighbor} \P\bigl( f_u \leq\eta_u
\leq M+1\mbox{,  for all } u\in N_{v_k} \mid \{\eta_{v_i} =
x_i\}_{0\leq i\leq k} \bigr) \geq10^{-(\Delta-1)}.
\end{equation}
Note that we do not condition on $\eta_{v_k}$ in the definition of
$f_w$, but we do condition on $\eta_{v_k}$ in
(\ref{eq-control-neighbor}). In order to prove
(\ref{eq-control-neighbor}), it suffices to show that for any $ f_w
\leq x_w \leq M + 1$ where $w\in B \subseteq N_{v_k}$, we have
%
\begin{eqnarray}
\label{eq-neighbor} \P\bigl( f_u \leq\eta_u \leq M+1
\mid\{\eta_{v_i} = x_i\}_{0\leq i\leq k} \cap\{
\eta_w = x_w\}_{w\in B}\bigr) \geq
\tfrac{1}{10}\nonumber\\[-8pt]\\[-8pt]
&&\eqntext{\mbox{for any } u\in N_{v_k} \setminus B.}
\end{eqnarray}
Write
$b_{y, z} = \P_y(\tau_{V_{k} \cup B} = \tau_z)$ for $y, z\in V$.
Applying Lemma \ref{lem-DGFF} and using the towering property of
conditional expectation, we get that
%
\begin{eqnarray}
\label{eq-towering-1} f_u &=& \E\bigl(\eta_u \mid\{
\eta_{v}\}_{v\in V_{k-1}}\bigr)= \E\bigl(\E\bigl(
\eta_u \mid\{\eta_v\}_{v\in V_{k} \cup B} \bigr) \mid\{
\eta_v\}_{v\in
V_{k-1}} \bigr)
\nonumber
\\
&=& \sum_{w\in V_{k} \cup B} b_{u, w} \E\bigl(
\eta_w \mid\{\eta_v\}_{v\in V_{k-1}}\bigr)
\\
&=& \sum_{w\in V_{k-1}} b_{u, w}
\eta_w + \sum_{w\in B \cup
\{v_k\} \setminus V_{k-1}} b_{u, w} \E
\bigl(\eta_w \mid\{\eta_v\}_{v\in V_{k-1}}\bigr)\nonumber
\end{eqnarray}
as well as that
%
\begin{equation}
\label{eq-towering-2} \E\bigl(\eta_u \mid\{\eta_{v}
\}_{v\in V_{k} \cup B}\bigr) = \sum_{w \in V_{k}
\cup B}
b_{u, w} \eta_w.
\end{equation}
Combined with (\ref{eq-x-k}), it follows that for any $\{x_w\}_{w\in
B}$, satisfying that $f_w \leq x_w\leq M+1$, for all $w\in B$, we have
\[
f_u \leq\E\bigl( \eta_u \mid\{\eta_{v_i} =
x_i\}_{0\leq i\leq k} \cap\{ \eta_w = x_w
\}_{w\in
B}\bigr) \leq M+1.
\]
Now an application of Lemma \ref{lem-DGFF} gives that
\[
\var\bigl(\eta_u \mid\{\eta_v\}_{v\in V_{k} \cup B}\bigr) =
R_{\mathrm
{eff}}(u, V_{k} \cup B) \leq R_{\mathrm{eff}}(u,
v_k) \leq1,
\]
where we used the fact that $u\sim v_k$. Recalling that $f_u \leq
M$, we can give a formal proof of (\ref{eq-neighbor}), and thereby establish
(\ref{eq-control-neighbor}) by a simple recursion.

A similar manipulation to (\ref{eq-towering-1}) using the towering
property of conditional expectation yields that
\begin{eqnarray*}
&&
\E\bigl(\eta_{v_k} \mid\{\eta_{v_i} = x_i
\}_{0\leq i <
k}\bigr) \\
&&\qquad = \E\bigl(\E\bigl(\eta_{v_k} \mid\{
\eta_{v_i} = x_i\}_{0\leq i < k}, \{ \eta_v
\}_{v\in N_{v_k}}\bigr) \mid\{\eta_{v_i} = x_i
\}_{0\leq i < k}\bigr)
\\
&&\qquad= \frac{1}{|N_{v_k}|} \sum_{w\in N_{v_k}} \E\bigl(
\eta_w \mid\{\eta_{v_i} = x_i
\}_{0\leq i < k}\bigr) = \frac{1}{|N_{v_k}|} \sum_{w\in
N_{v_k}}
f_w = \sum_{i<k} a_{k, i}
x_i.
\end{eqnarray*}
Recalling that $f_w \leq M$ for all $w\in V$, we have
\[
\sum_{i<k}a_{k, i} |M - x_i|
= \sum_{i<k}a_{k, i} (M - x_i)
= \frac
{1}{|N_{v_k}|} \sum_{w\in N_{v_k}}(M -
f_w).
\]
Take $\{x_w\}_{w\in N_{v_k}}$ such that $f_w\leq x_w \leq M+1$ for all
$w\in N_{v_k}$. Since $f_w\leq M$ for all $w\in N_{v_k}$, we have
\[
\frac{1}{|N_{v_k}|}\sum_{w\in N_{v_k}} |M - x_w|
\leq1 + \frac
{1}{|N_{v_k}|} \sum_{w\in N_{v_k}}(M -
f_w) = 1 + \sum_{i<k}a_{k, i}
|M - x_i|.
\]
By (\ref{eq-x-k}) and the fact that $ 1\wedge\frac{1}{x} \leq1
\wedge\frac{2}{x+1}$ for all $x > 0$, we obtain that
\[
M \leq x_k \leq M + 2 \epsilon\biggl(1 \wedge\frac{|N_{v_k}|}{\sum
_{w\in N_{v_k}}|M-
x_w|}
\biggr) \leq M + 2 \epsilon\biggl(1 \wedge\frac{\Delta}{\sum_{w\in
N_{v_k}}|M- x_w|} \biggr).
\]
Combined
with (\ref{eq-control-v-k}) and (\ref{eq-control-neighbor}), it
follows that for any $(x_0,\ldots, x_{k-1}) \in(-\infty, M)^k$,
\begin{eqnarray*}
&&
\P \biggl(M \leq\eta_{v_{k}} \leq M + 2\epsilon\biggl(1 \wedge
\frac{\Delta}{\sum_{u\in N_{v_k}}
|M-\eta_{u}|} \biggr) \Bigm|\{\eta_{v_i} = x_i
\}_{0\leq i<k} \biggr)
\\
&&\qquad\geq\frac{\epsilon}{5} \cdot10^{-(\Delta-1)}\cdot\P(\eta_{v_{k}} \geq
M \mid\eta_{v_i} = x_i \mbox{ for } 0\leq i < k).
\end{eqnarray*}
Combined with (\ref{eq-decomoposition-1}) and
(\ref{eq-decomoposition-2}), it follows that
\[
\P(\Xi) \geq\frac{\epsilon}{5} \cdot10^{-(\Delta-1)} \cdot\P\Bigl(\sup
_v \eta_v \geq M \Bigr),
\]
completing the proof.
\end{pf}

In the particular case that $M$ is the median of $\sup_v \eta_v$, the
preceding proposition gives that
%
\begin{equation}
\label{eq-gff-detection} \P\biggl(\exists v\in V\dvtx  M \leq\eta_{v}
\leq M + \epsilon\biggl(1 \wedge\frac{\Delta}{\sum_{u\in N_v}|M-\eta
_{u}|} \biggr) \biggr) \geq
\frac{\epsilon}{ 10^{\Delta}}.
\end{equation}

We remark that both (\ref{eq-GFF-exp}) and (\ref{eq-GFF-var}) in the
sequential decomposition of GFF are of crucial importance to our
proof: (\ref{eq-GFF-exp}) ensures that the conditional mean of any
variable that is being revealed, is less than the target value $M$
as long as none of the previous values exceeds $M$;
(\ref{eq-GFF-var}) guarantees that we can reveal the GFF in a way
such that the conditional standard deviation of each variable is
bounded by 1.

\section{Reconstructing random walks from local
times}\label{secreconstruction}

In this section, we study the reconstruction of random walks from
local times. Roughly speaking, conditioning on the local times, the
embedded discrete-time random walk should be biased to those paths
that are more likely to fulfill the desired local times. The goal is
to understand such bias implied in the local times.

In Section \ref{secreconstruction-set-up}, we set up the general
framework for the study of the conditioned measure of embedded walks
given local times and demonstrate a connection to enumeration of
Eulerian circuits. In Section \ref{secthinpoints}, we focus on
the number of visits to a certain vertex, and give an upper bound
assuming a small local time at this vertex as well as its neighbors.
Section \ref{secsprinkling} contains a sprinkling argument
which, together with results obtained in
Section \ref{secdetection}, proves the Theorem \ref{thm-tau-cov}.

\subsection{Random walks, local times and Eulerian
circuits}\label{secreconstruction-set-up}

Let $G=G(V)$ be a network with conductance $c_{u, v}$ on edge $(u,
v)$. Let $X_t$ be a continuous-time random walk on $G$ associated with
$\{c_{u, v}\}$. Define a
sequence of stopping times $\tau_k$ in the following way:
\[
\tau_0 = 0,\qquad \tau_{k} = \inf\{t>\tau_{k-1}\dvtx
X_t \neq X_{\tau_{k-1}}\} \qquad\mbox{for } k\geq1.
\]
Define $K = \max\{k\dvtx  \tau_k \leq\tau(t)\}$. We consider the
embedded discrete-time random walk $(S_k)$, which is defined to be
$S_k = X_{\tau_k}$ for $k=0, 1,\ldots, K$. Note that $S_0 = S_K =
v_0$. Let $\mathcal{P}$ be the path of random walk $(S_k)$ up to
time $K$, and let $\Omega$ be the space of paths which start and end
at $v_0$ and visit every vertex in the graph. For $P \in\Omega$ and
$u\neq v$,
define $k_{u, v} = k_{u, v}(P)$ to be the number of times
that path $P$ traverses the directed edge $\langle u, v \rangle$, and
define $k_v = k_v(P) = \sum_{u\sim v} k_{u, v}(P)$ to be the number of
times that path $P$ visits
$v$.

The central task of this section is to reconstruct path
$\mathcal{P}$ generated by random walk $(S_k)$ conditioning on the
event that
%
\begin{equation}
\label{eq-def-Gamma} \Gamma= \bigl\{L^v_{\tau(t)} =
\ell_v: \mbox{ for } v\in V \setminus\{v_0\}\bigr\}.
\end{equation}
For convenience of notation, we
write $t = \ell_{v_0} $.

\begin{lemma}\label{lem-random-walk-path}
Write $\check{c}_{v} = c_v - c_{v, v}$ for all $v\in V$. We have that
\[
\mu(\mathcal{P} = P, \Gamma) =\mathrm{e}^{-\check{c}_{v_0}t} \frac
{t^{k_{v_0}}}{k_{v_0} !} \cdot
\prod_{u\neq v}c_{u, v}^{k_{u, v}} \cdot\prod
_{v\neq v_0}\frac{\ell_v^{k_v-1}
\mathrm{e}^{-\check{c}_v\ell_v}}{ (k_v-1)!}.
\]
Here $\mu(\mathcal{P} = P, \Gamma) = \P(\mathcal{P} = P) \mu
(\Gamma\mid\mathcal{P} = P)$ with $\mu(\Gamma\mid\mathcal{P} =
P)$ being the conditional density of the local times with respect to
the Lebesgue measure, given that $\mathcal{P} = P$.
\end{lemma}

\begin{pf}
It is clear that conditioning on $\{\mathcal{P} = P\}$, we have that
for all $v\neq v_0$,
\[
\bigl(L^v_{\tau(t)} \mid\mathcal{P} = P\bigr) \stackrel{\mathrm{law}}
{=} \frac{1}{c_v}\sum_{i=1}^{k_v}
Y_{v, i},
\]
where $\{Y_{v, i}\}_{v\in V, i\in\N}$ is a collection of
independent exponential variables with $\E Y_{v, i} =
\frac{c_v}{\check{c}_{v}}$. Therefore, we have for all $v\neq v_0$,
\[
\bigl(L^v_{\tau(t)} \mid\mathcal{P} = P\bigr) \stackrel{\mathrm{law}}
{=} \frac{1}{\check{c}_v}\sum_{i=1}^{k_v}
Z_{v, i},
\]
where $\{Z_{v, i}\}_{v\in V, i\in\N}$ is a collection of i.i.d.
standard exponential variables. This implies that
%
\begin{eqnarray}
\label{eq-local-times}
\mu( \Gamma\mid\mathcal{P} = P) &=& \mu\bigl
(L^v_{\tau(t)}
= \ell_v\mbox{,  for } v\in V\setminus\{v_0\} \mid
\mathcal{P} = P\bigr) \nonumber\\[-8pt]\\[-8pt]
&=& \prod_{v\neq v_0} \bigl(
\check{c}_v g(k_v, \check{c}_v
\ell_v) \bigr),\nonumber
\end{eqnarray}
where $g(k, x) =
\frac{x^{k-1} \mathrm{e}^{-x}}{(k-1)!}$ is the density at $x$ of a
Gamma variable with parameter $(k, 1)$, and the factor $\check{c}_v$
before $g(\cdot, \cdot)$ comes from change of variables. In
addition, by definition of continuous-time random walks, we see that
the number of excursions $K_{v_0}$ at $v_0$ accumulated up to
$\tau(t)$ is distributed as a Poisson variable with mean
$\check{c}_{v_0} t$, and it is independent of the realization of the
excursions. Therefore,
%
\begin{equation}
\label{eq-path} \P(\mathcal{P} = P)= \P(K_{v_0} = k_{v_0})
\cdot\prod_{u \neq v}p_{u, v}^{k_{u, v}} =
\mathrm{e}^{-\check{c}_{v_0}t} \frac{(\check
{c}_{v_0}t)^{k_{v_0}}}{k_{v_0}!} \cdot\prod
_{u\neq v}p_{u,
v}^{k_{u, v}},
\end{equation}
where $\prod_{u\neq v}p_{u, v}^{k_{u, v}}$ counts the probability
for the embedded random walk to follow path $P$. Combining
(\ref{eq-local-times}) and (\ref{eq-path}), we conclude that
\[
\mu(\mathcal{P} = P, \Gamma) = \mathrm{e}^{-\check{c}_{v_0}t} \frac
{t^{k_{v_0}}}{k_{v_0} !}
\cdot\prod_{u\neq v}c_{u, v}^{k_{u, v}} \cdot
\prod_{v\neq
v_0}\frac{\ell_v^{k_v-1} \mathrm{e}^{-\check{c}_v\ell_v}
}{(k_v-1)!},
\]
completing the proof.
\end{pf}
In light of Lemma \ref{lem-random-walk-path}, we could write
\[
\mu(\mathcal{P} = P, \Gamma) = \biggl(\mathrm{e}^{-\check
{c}_{v_0}t} \prod
_{u\neq v_0} \frac{\mathrm{e}^{-\check{c}_v\ell
_v}}{\ell_v} \biggr)\cdot\prod
_{u\neq v}c_{u, v}^{k_{u, v}} \cdot\frac{t^{k_{v_0}}}{k_{v_0} !}
\cdot\prod_{v\neq v_0}\frac{\ell_v^{k_v}
}{ (k_v-1)!}.
\]
Since the prefactor in the parentheses above is common for every path
$P$ (assuming $\Gamma$ is fixed), the following corollary is now immediate.
%
\begin{cor}\label{cor-random-walk-path}
For real numbers $\ell_v \geq0$ for $v\in V$, let $\Gamma$ be
defined as in~(\ref{eq-def-Gamma}). Write $t = \ell_{v_0}$. For any
$P\in\Omega$, write
%
\begin{equation}
\label{eq-W-P} W_P = \prod_{u\neq v}c_{u, v}^{k_{u, v}}
\frac{t^{k_{v_0}}}{k_{v_0} !} \prod_{v\neq v_0}\frac{\ell
_v^{k_v}}{(k_v-1)!}.
\end{equation}
Then, for all $P \in\Omega$,
\[
\P(\mathcal{P} = P \mid\Gamma) = \frac{W_P}{Z},
\]
where $Z$ is a
normalizing constant depending on $\Gamma$.
\end{cor}

The next fact follows immediately from
Corollary \ref{cor-random-walk-path} and reversibility of random
walks.
%
\begin{claim}\label{claim-reverse-cycle}
Let $P \in\Omega$ be a random walk path, and suppose $P$ consists
of three parts, $P_1, C, P_2$, in an order where $C$ is a cycle and $P_1,
P_2$ are paths. Let $\overleftarrow{C}$ be the reversed cycle of
$C$, and let $\tilde{P}$ be a path consisting of $P_1,
\overleftarrow{C}, P_2$, in order. Then $W_P = W_{\tilde{P}}$.
\end{claim}

We now explore a connection between random walk paths and Eulerian
graphs. Given a collection $\{j_{u, v}\dvtx  u, v\in V\}$, we let
$\mathcal{G} = \mathcal{G}(V)$ be a multiple directed graph where
the multiplicity of directed edge $\langle u, v\rangle$ is $j_{u,
v}$. We say $\mathcal{G}$ is Eulerian if there is a Eulerian circuit
for the graph $\mathcal{G}$, that is, a circuit which traverses every directed
edge in the graph exactly once. A~classical argument says that a
directed graph $\mathcal{G}$ is Eulerian if and only if it is
connected, and the in-degree is equal to the out-degree for every
vertex in $\mathcal{G}$. Clearly, if $j_{u, v} = k_{u, v}(P)$ for a
certain $P\in\Omega$, the associated graph $\mathcal{G}$ is
Eulerian. Denote by $\operatorname{ec}(\mathcal{G})$ the number of
Eulerian circuits for graph~$\mathcal{G}$ (where the circuits that
are identical up to cyclic translations are counted only once). For
$v\in\mathcal{G}$, let $\operatorname{ec}_v(\mathcal{G})$ be the number
of Eulerian circuits started at~$v$ (where we do distinguish
circuits obtained from cyclic translations). Note that
$\operatorname{ec}_v(\mathcal{G}) = \mathrm{deg}_v \cdot
\operatorname{ec}(\mathcal{G})$, where $\mathrm{deg}_v$ is the in-degree
(equivalently, out-degree) of vertex $v$.

\begin{claim}\label{claim-traverse}
Let $\mathcal{G}$ be a multiple directed graph associated with $\{
j_{u, v}\}_{u, v\in
V}$ such that $j_{v, v} = 0$ for all $v\in V$, and suppose that
$\mathcal{G}$ is Eulerian. Define
\[
\Omega(\mathcal{G}) = \Omega\bigl(\{j_{u, v}\}_{u, v\in V}\bigr) =
\bigl\{P\in\Omega\dvtx  k_{u, v}(P) = j_{u, v}\mbox{,  for } u, v \in
V\bigr\}.
\]
Then $|\Omega(\mathcal{G})| = \operatorname{ec}_{v_0}(\mathcal{G}) \cdot
(\prod_{u, v}j_{u, v}!)^{-1}$.
\end{claim}

\begin{pf}
Consider the multiple directed graph $\mathcal{G}$. We see that each
Eulerian circuit started at $v_0$ on $\mathcal{G}$ induces a
legitimate
path $P\in\Omega$, and a path $P\in
\Omega$ corresponds to a number of Eulerian circuits. Since the
multiple edges in $\mathcal{G}$ are not distinguishable in the
path, the mapping from Eulerian circuits to $\Omega(\mathcal{G})$ has
multiplicity $\prod_{u, v\in V} j_{u, v}!$. This completes
the proof.
\end{pf}

We have the following classic result on the enumeration of Eulerian
circuits for directed Eulerian graphs, known as BEST theorem that
was originally proved by Aardenne-Ehrenfest and de Bruijn
\cite{BE51} as a variation of an earlier result of Smith and Tutte
\cite{TS41}.

\begin{theorem}[(\cite{BE51,TS41})]\label{thm-BEST}
Let $\mathcal{G} = (V, E)$ be a multiple directed Eulerian graph. Then
for any $w\in V$
\[
\operatorname{ec}(\mathcal{G}) = \operatorname{ar}_w(\mathcal{G}) \prod
_{v\in
V} (\mathrm{deg}_v - 1)!,
\]
where $\operatorname{ar}_w(\mathcal{G})$ is the number of arborescences,
which are
directed trees such that there exists a unique path towards the
vertex $w$ for every $v\in V$ and $v \neq w$.
\end{theorem}
\begin{remark*}
It is implied from the preceding theorem that
for Eulerian graphs,
%
\begin{equation}
\label{eq-referee} \operatorname{ar}_w(\mathcal{G}) =
\operatorname{ar}_v(\mathcal{G})\qquad\mbox{for all }w, v\in\mathcal{G}.
\end{equation}
\end{remark*}

The next corollary is an immediate consequence of Corollary \ref
{cor-random-walk-path}, Claim \ref{claim-traverse} and Theorem \ref{thm-BEST}.

\begin{cor}\label{cor-random-eulerian}
Let $\mathcal{G}$ be a Eulerian graph associated with $\{j_{u,
v}\}_{u, v\in V}$. Then
%
\begin{equation}
\label{eq-random-eulerian} \P\bigl(\mathcal{P} \in\Omega(\mathcal{G})
\mid\Gamma
\bigr) = \frac{\operatorname{ar}_{v_0}(\mathcal{G})}{Z} \prod_{u \neq v}
\frac{(\sqrt{\ell_u\ell_v} c_{u, v})^{j_{u, v}}}{j_{u, v}
!},
\end{equation}
where $Z$ is a normalizing constant.
\end{cor}

The next claim on the enumeration of arborescences will be useful.

\begin{claim}\label{claim-enumeration-arborescence}
Let $\mathcal{G}$ and $\mathcal{G}'$ be two directed Eulerian graphs
over vertex set $V \ni v_0$. Denote by $j_{u, v}$ and $j'_{u, v}$
the multiplicity of edge $\langle u, v\rangle$ in graph
$\mathcal{G}$ and $\mathcal{G}'$, respectively. Suppose that $j'_{u,
v} \geq1$ for all $u, v\in V$ if and only if $j_{u, v} \geq1$. Assume also
that $j_{u, v} \geq j'_{u, v}$ for all $u, v\in V$. Then
\[
\frac{\operatorname{ar}_{v_0} (\mathcal{G}')}{\mathrm{ar}_{v_0}(\mathcal
{G})} \geq\prod_{u\neq v: j_{u, v}\geq1}
\frac{j'_{u, v}}{j_{u,
v}}.
\]
\end{claim}
\begin{pf}
Consider an arborescence $T$ for the complete directed graph on vertex
set $V$. Denote by $\mathcal{A}_T$ and $\mathcal{A}'_T$ the set of
arboresences that correspond to $T$ in $\mathcal{G}$
and~$\mathcal{G}'$, respectively. By correspondence, we mean that they
are the same if we identify all the multiple edges. That is, we say $T$
corresponds to $T'$ if and only if for every edge $e\in T$, there
exists an edge $e' \in T'$ where $e$ and $e'$ share the same starting
and ending points. Since $j'_{u, v} \geq1$ whenever $j_{u, v} \geq1$ by
our assumption, we see that $|\mathcal{A}'_T|\geq1$ as long as
$|\mathcal{A}_T| \geq1$. Furthermore, it is clear that the cardinality
of $\mathcal{A}_T$ ($\mathcal{A}'_T$) is determined by the multiplicity
of the edges that appear in $T$. More precisely,
\[
|\mathcal{A}_T| = \prod_{\langle u, v \rangle\in T}
j_{u, v},
\]
and a similar equality holds for $\mathcal{A}'_T$. Therefore, for every $T$,
we have
\[
\frac{|\mathcal{A}'_T|}{ |\mathcal{A}_T|} \geq\prod_{u\neq v}
\frac{j'_{u, v}}{j_{u, v}}.
\]
Summing over arborescence
$T$, we can then deduce the claim.
\end{pf}

\subsection{Thin points of random walks}\label{secthinpoints}

In this subsection, we study the probability distribution on
$\Omega$ conditioning on $\Gamma$. For $v \in V$, we call $v$
an $m$-thin point of random walk path $P\in\Omega$ if $k_v(P) \leq
m$. We wish to lower bound the probability that the random walk path
has $v$ as an $m$-thin point for suitable $m\in\N$, conditioned on
the local times.\vadjust{\goodbreak}

We first demonstrate that the number of traverses over edge $\langle
u, v \rangle$ cannot be too different from the number of traverses
over edge $\langle v, u\rangle$. Note that we use $\langle u, v
\rangle$ to denote directed edges.

\begin{lemma}\label{lem-traverse-balance}
Let $\nu$ be the probability measure conditioning on the event $\Gamma
$. For all $u\neq v$ and $k\geq184$, we have
\begin{eqnarray*}
&&\nu\bigl(k_{u, v}(P) + k_{v, u}(P) = k, \bigl|k_{u, v}(P) -
k_{v, u}(P)\bigr| \geq k/3\bigr) \\
&&\qquad\leq\tfrac{1}{2} \nu
\bigl(k_{u, v}(P) + k_{v, u}(P) = k\bigr).
\end{eqnarray*}
\end{lemma}
\begin{pf}
Fix arbitrary $u\neq v$ and $k\geq184$. Consider $P\in\Omega$
with $k_{u, v}(P) + k_{v, u}(P) = k$. We can decompose $P$ into a
sequence of $(P_1, C_1,\ldots, C_\ell, P_2)$ for a certain $\ell=
\ell(P) \in\N$, where $C_i$ is a cycle containing either one edge
$\langle u, v\rangle$ or one edge $\langle v, u\rangle$ or one pair
of them, and $P_1, P_2$ are paths such that $(P_1, P_2)$ forms a
cycle containing at most one traverse between $u$ and $v$ (note that a
cycle in this paper is simply a path such that the starting and ending
points are the same, with no self-avoiding constraints posed). Indeed,
we select the following specified manner for the decomposition:
along path $P$, let $P_1$ be the segment of $P$ from $v_0$ to the
first encounter of a vertex $w\in\{u, v\}$; we
then continue searching along path $P$ and let $C_1$ be the segment
of $P$ until for the first time the random walk goes back to $w$
after experiencing a traverse between $u$ and $v$ (hence $C_1$ is a
cycle); we repeat this procedure to obtain $C_i$ until no such
cycles exist. The last segment of path leading back to $v_0$ is
then defined to be $P_2$. Note that $(k-1)/2\leq\ell\leq k$.

For $1\leq i\leq\ell$, assume that $C_i$ contains cycles $C_{i, 1},\ldots, C_{i, \ell_i}$ in this order where $C_{i, j} = w, x^{(i,
j)}_{1},\ldots, x^{(i, j)}_{r_{i, j} + 1} = w$ and $x_{m}^{(i, j)}
\neq w$ for all $1\leq m\leq r_{i, j}$. Define the reverse of $C_i$ to
be the cycle consisting of $\overleftarrow{C}_{i, 1},\ldots,
\overleftarrow{C}_{i, \ell_i}$ in this order where $\overleftarrow
{C}_{i, j} = w, x^{(i, j)}_{r_{i,j}},\ldots, x^{(i, j)}_{1}, w$. Let
$\ell'\leq\ell$ be the number of cycles in $\{C_1,\ldots, C_\ell\}
$ such that the reverse is different from the original one, and we
assume that these cycles are $C'_1,\ldots, C'_{\ell'}$. We write
$P'\sim P$ if $P' \in\Omega$ can be obtained from $P$ by
reversing a subset of the cycles $\{C'_1,\ldots, C'_{\ell'}\}$, as in
the statement of Claim \ref{claim-reverse-cycle}. It is clear that
$\sim$ is an equivalent relation on $\Omega$ and hence generates a
partition. Denote by $\Omega_P = \{P'\in\Omega\dvtx  P' \sim P\}$. We
see that $|\Omega_P| = 2^{\ell'}$, and $k_{v, u}(P') + k_{u, v}(P') =
k$ for $P'\in\Omega_P$. By Claim \ref{claim-reverse-cycle},
we have
\[
\nu(P) = \nu\bigl(P'\bigr) \qquad\mbox{for all }P'\sim P.
\]
Let $\chi_i =
\one_{\langle u, v\rangle\in C'_i}$ and $\chi'_i = \one_{\langle v,
u\rangle\in C'_i}$. Denote by $D_i \in\{-1, 1\}$ the direction of
cycle $C'_i$ (we use the convention that all $D_i = 1$ for $P$). A
simple application of Chernoff bound gives that
\[
\biggl| \biggl\{(D_1,\ldots, D_{\ell'}) \in\{-1, 1
\}^{\ell'}\dvtx  \biggl|\sum_{i} \bigl(
\chi_i - \chi'_i\bigr) D_i \biggr|
\geq\frac{k}{4} \biggr\} \biggr| \leq2^{\ell'-1}.
\]
This implies that
\[
\nu\bigl(\bigl\{P'\in\Omega_P\dvtx  \bigl|k_{u, v}
\bigl(P'\bigr) - k_{v, u}\bigl(P'\bigr)\bigr| \geq
k/3\bigr\}\bigr) \leq\nu(\Omega_P)/2.
\]
The proof is completed by summing over all classes $\Omega_P$ for
$k_{u, v}(P) +\break  k_{v, u}(P) = k$.
\end{pf}

We next turn to analyze the number of traverses between vertices $u$
and~$v$.

\begin{lemma}\label{lem-traverses-edge}
Fix $u\neq v$ and suppose that $\ell_u\ell_v c^2_{u, v} \leq
1/16$. Then for any \mbox{$k \geq184$}, we have
\[
\nu\bigl(k_{u, v}(P) + k_{v, u}(P) = k+1\bigr) \leq\nu
\bigl(k_{u, v}(P) + k_{v,
u}(P) = k-1\bigr)/2.
\]
\end{lemma}
\begin{pf}
By Lemma \ref{lem-traverse-balance}, we have that
%
\begin{eqnarray}
\label{eq-balance-traverse}
&& \nu\bigl(k_{u, v}(P) + k_{v, u}(P) =
k+1, \bigl|k_{u, v}(P) - k_{v, u}(P)\bigr| \geq(k+1)/3\bigr)\nonumber\\[-8pt]\\[-8pt]
&&\qquad \leq\nu
\bigl(k_{u, v}(P) + k_{v, u}(P) = k+1\bigr)/2.\nonumber
\end{eqnarray}
Therefore, it
suffices to bound the measure for the set
\[
\Omega_{k+1} \deq\bigl\{P\in\Omega\dvtx  k_{u, v}(P) +
k_{v, u}(P) = k+1, \bigl|k_{u, v}(P) - k_{v, u}(P)\bigr|
\leq(k+1)/3\bigr\}.
\]
For $P\in\Omega$, denote by $\mathcal{G}_P$ the directed Eulerian
graph generated by $P$ [i.e., the multiplicity of edge $\langle x,
y\rangle$ in $\mathcal{G}$ is $k_{x, y}(P)$ for $x, y\in V$]. For
$P, P'\in\Omega$, we say $P \sim_{\mathcal{G}} P'$ if
$\mathcal{G}_P = \mathcal{G}_{P'}$. Clearly, $\sim_{\mathcal{G}}$ is
an equivalent relation and hence generates a partition over
$\Omega_{k+1}$. Now take $P\in\Omega_{k+1}$ and denote by
$\mathcal{G}$ the generated Eulerian graph. We next study the
equivalent class of $P$.

By definition of $\Omega_{k+1}$ and the fact that $k\geq12$, we see
that $k_{u, v}(P),\break  k_{v, u}(P) \geq2$. Hence, we can obtain an
Eulerian graph $\mathcal{G}'$ from $\mathcal{G}$ by reducing the
multiplicity of both $\langle u, v \rangle$ and $\langle v, u
\rangle$ by 1. By Claim \ref{claim-enumeration-arborescence}, we get
that $\operatorname{ar}_{v_0}(\mathcal{G}) \leq4
\operatorname{ar}_{v_0}(\mathcal{G}')$. Now applying
Corollary \ref{cor-random-eulerian} and using our assumption on $\ell
_u\ell_v$, we get that
$\nu(\Omega(\mathcal{G})) \leq\nu(\Omega(\mathcal{G}'))/4$. Summing
over all $\mathcal{G}$ generated by paths in $\Omega_{k+1}$, we get
that
\[
\nu(\Omega_{k+1}) \leq\nu\bigl(k_{u, v}(P) +
k_{v, u}(P) = k-1\bigr)/4.
\]
Combined with (\ref{eq-balance-traverse}), the required inequality
follows.
\end{pf}

We then arrive at the following consequence.

\begin{prop}\label{prop-thin-path}
Consider a network $G = G(V)$ with $v_0\in V$ a fixed $v\in V$. Let
$N_v =
\{u\neq v\dvtx  c_{u, v} > 0\}$. Suppose that $\{\ell_w\}_{w\in V}$ are
positive numbers such that $\ell_u\ell_v c^2_{u, v} \leq1/16$ for
all $u\in N_v$ (otherwise we already have a unvisited vertex). Define
$\Gamma= \{L^w_{\tau(t)} = \ell_w$ for
all $w\in V\}$. Denoting by $\mathcal{P}$ a random path for the
embedded walk up to time $\tau(t)$, we have
\[
\P\bigl(k_v(\mathcal{P}) \leq1118|N_v| \mid\Gamma\bigr)
\geq1/2.
\]
\end{prop}
\begin{pf}
An application of Lemma \ref{lem-traverses-edge} yields that
$\E(k_{u, v}(\mathcal{P})\mid\Gamma) \leq559$ for all $u\in N_v$,
and thereby $\E(k_v(\mathcal{P}) \mid\Gamma) \leq559 |N_v|$. The
proof is completed by a simple application of Markov inequality.
\end{pf}

\begin{remark*}
The bounded-degree assumption was made in order
for the preceding proposition to be useful---the bound that was proved
on the number of visits to $v$ grows linearly with the degree of $v$,
and thus stopped being useful if the degree of $v$ is unbounded in the
sequence of graphs. We were hoping that a more careful reconstruction
argument could yield an upper bound that is independent of the degree,
but it seems that this could not be achieved by the current method
which considers the number of traverses from all the neighboring edges
separately.
\end{remark*}

\subsection{A sprinkling argument}\label{secsprinkling}

In this subsection, we establish Theorem \ref{thm-tau-cov} based on
results developed in previous sections. We first demonstrate the
existence of thin points for random walks with nonnegligible
probability. For a continuous-time random walk $(X_s)$,
we denote by $K_v(s)$ the number of visits to vertex $v\in V$ up to
time $s$, for the corresponding embedded discrete-time random walk.
That is to say, $K_v(s)$ is the maximal number $k$ such that there
exists $s_0 < s_1 < s_2 < \cdots< s_{2k} \leq s$ with $X_{s_i} = v$ for
even $i$ and $X_{s_i} \neq v$ for odd $i$.

\begin{prop}\label{prop-thin-points}
For a graph $G = G(V, E)$ with maximal degree bound\-ed by $\Delta$,
let $\{\eta_v\}_{v\in V}$ be a GFF on $G$ with $\eta_{v_0} = 0$ for
a certain $v_0 \in V$. Denote by $M$ the median of $\sup_v \eta_v$,
and write $t = M^2/2$. Let $\tau(t)$ be defined as in~(\ref{eqinverselt}). Then
\[
\P_{v_0}\bigl(\exists v\in V\dvtx  K_v\bigl(\tau(t)\bigr)
\leq1118 \Delta\bigr) \geq\frac
{1}{8\Delta\cdot10^\Delta}.
\]
\end{prop}
\begin{pf}
Applying (\ref{eq-gff-detection}) with $\epsilon=
\frac{1}{4\Delta}$, we obtain that
\[
\P\biggl(\exists v\in V\dvtx  M \leq\eta_{v} \leq M + \frac{1}{4\Delta}
\biggl(1 \wedge\frac{\Delta}{\sum_{u\in N_v}|M-\eta_{u}|} \biggr) \biggr
) \geq\frac{1}{4 \Delta\cdot10^{\Delta}}.
\]
In particular, this implies that with probability at least
$\frac{1}{4 \Delta\cdot10^{\Delta}}$, there exists $v\in V$ such
that
\[
|\eta_v - M| \cdot|\eta_u - M| \leq1/4\qquad \mbox{for all
} u\in N_v.
\]
In view of Theorem \ref{thmrayknight}, we see that
\[
\bigl\{L^v_{\tau(t)}\dvtx  v\in V\bigr\} \preceq\bigl\{(
\eta_v - \sqrt{2t})^2\dvtx  v\in V\bigr\}.
\]
Altogether, we obtain that with probability at least
$\frac{1}{4 \Delta\cdot10^{\Delta}}$, there exists $v\in V$ such
that
\[
L^v_{\tau(t)} L^u_{\tau(t)} \leq1/16\qquad
\mbox{for all } u\in N_v.
\]
Combined with Proposition \ref{prop-thin-path}, the desired lower
bound on the probability follows.
\end{pf}

At this point, we employ a sprinkling argument to complete the
proof. The basic intuition is that there should be a nonnegligible
chance that the random walk fails to cover the graph at time
$\tau((1-\epsilon)t)$, given that the random walk barely covers the
graph at time $\tau(t)$.

\begin{pf*}{Proof of Theorem \ref{thm-tau-cov}}
Let $t = M^2/2$. Denote by $F$ the event $\{\exists v\in V\dvtx
K_v(\tau(t)) \leq1118 \Delta\}$. Proposition \ref{prop-thin-points}
asserts that $\P_{v_0}(F) \geq1/(8\Delta\cdot10^\Delta)$. We
next condition on the event $F$.

Let $\mathcal{C} = \{C_1,\ldots, C_N\}$ be the (multiple) set of
excursions at the origin $v_0$, where $N$ is the number of total
excursions occurring at $v_0$ up to time $\tau(t)$. We emphasize
that we define $\mathcal{C}$ to be the set of excursions, without
distinguishing the orderings among the excursions. In particular,
$C_i$ is not necessarily the $i$th excursion that occurs in the
random walk. Furthermore, the ordering of the occurrences of these
excursions forms a uniformly random permutation. Denote by $0\leq
T_i\leq t$ the local time at $v_0$ when the excursion $C_i$ occurs.
A crucial observation is that, conditioning on $\mathcal{C}$ as well
as event $F$, the random times $T_i$ are i.i.d. uniformly
distributed over $[0, t]$, since $\{T_i\}$ arises from a Poisson
point process on $[0,t]$.

Recall the definition of $F$, we can now select a vertex $v\in V$
such that there are at most $1118\Delta$ excursions that ever visited
$v$. Let $I = \{1\leq i\leq N\dvtx  v\in C_i\}$. It is now clear that
\[
\P\bigl(\forall i\in I\dvtx  T_i \geq(1 - \epsilon)t \mid F, \mathcal{C}
\bigr) = \epsilon^{|I|} \geq\epsilon^{1118\Delta}.
\]
This implies that
\[
\P\bigl(\tau_{\mathrm{cov}} \geq\tau\bigl((1-\epsilon)t\bigr) \mid
F\bigr) \geq
\epsilon^{1118\Delta}.
\]
Combined with the lower bound on $\P(F)$, it follows that
%
\begin{equation}
\label{eq-for-referee} \P\bigl(\tau_{\mathrm{cov}} \geq\tau\bigl
((1-\epsilon)t
\bigr)\bigr) \geq\frac
{\epsilon^{1118\Delta}}{8\Delta\cdot10^\Delta}.
\end{equation}
Now by Proposition \ref{prop-DLP} and assumption (\ref
{eq-assumption}) as well as the commute time identity (\ref
{eq-commute-time}), we see that
\[
\frac{R}{(\E\sup_v \eta_v)^2} \leq\frac{2t_{\mathrm
{hit}}}{t_{\mathrm{cov}}} \cdot(1 + C\sqrt{t_{\mathrm
{hit}}/t_{\mathrm{cov}}})
\leq\frac{4 \epsilon^4}{10^4 \Delta^2},
\]
where $R$ is the diameter in resistance metric. Note that by Lemma
\ref{lem-gaussian-concentration}, we have $M \geq(1- \epsilon/4) \E
\sup_v
\eta_v$. Thus
\[
2(1-2\epsilon) t |E| = (1-2\epsilon)|E|M^2 \geq(1-3\epsilon) |E|
\Bigl(\E\sup_v\eta_v\Bigr)^2.
\]
Therefore,
\begin{eqnarray*}
&&
\Bigl\{\tau\bigl((1-\epsilon)t\bigr) \leq(1-3\epsilon) |E| \Bigl(\E\sup
_v\eta_v\Bigr)^2\Bigr\}\\
&&\qquad \subseteq
\bigl\{\tau\bigl((1 - \epsilon)t\bigr) - 2(1 - \epsilon)t |E| \leq-2
\epsilon t
|E|\bigr\}.
\end{eqnarray*}
Applying Lemma \ref{lem-concentration-tau-t}, aiming at a
concentration of $\tau((1 - \epsilon)t)$ with a choice of $\lambda=
2\epsilon t/(\sqrt{tR} + R) \geq10 \Delta/\epsilon^2$ (where
$\lambda$ is a parameter in the statement of the Lemma
\ref{lem-concentration-tau-t}), we obtain that for $\epsilon\leq10^{-4}$
(note that $\Delta\geq2$)
\[
\P\Bigl(\tau\bigl((1-\epsilon)t\bigr) \leq(1-3\epsilon) |E| \Bigl(\E
\sup
_v\eta_v\Bigr)^2\Bigr) \leq6
\mathrm{e}^{-\lambda/16} \leq\frac{\epsilon
^{1118\Delta}}{32\Delta\cdot10^\Delta}.
\]
Combined with (\ref{eq-for-referee}), the conclusion of the theorem
follows with a choice of $\delta=
\frac{(\epsilon/3)^{1118\Delta}}{32\Delta\cdot10^\Delta}$.
\end{pf*}

\section{Discussions and future directions}\label{secdiscussion}

Our work (obviously) reinforces a number of questions on cover
times posed in \cite{DLP10}, including the asymptotics and
exponential concentration for cover times on general graphs. In what
follows, we discuss additional three questions motivated by the
current work.

\textit{Deterministic approximation scheme for Gaussian free
field.} The resolution of deterministic polynomial-time
$O(1)$-approximation for cover times on general graphs \cite{DLP10},
naturally raises the question of designing a deterministic
polynomial-time approximation scheme (DPTAS). That is, a
deterministic algorithm which takes $\epsilon>0$ as a parameter and
approximates the cover time up to a factor of $1+\epsilon$ in
polynomial-time (where the power of the polynomial depends on
$\epsilon$). This question was solved for general trees \cite{FZ09}
using dynamic programming. Our work confirms that cover times can be
recovered from GFF with a precision up to $1+o(1)$ for general trees
and bounded degree graphs, assuming $t_{\mathrm{hit}} =
o(t_{\mathrm{cov}})$. Therefore, a DPTAS for the supremum of GFF
will immediately give a DPTAS for cover times on bounded degree
graphs, and plausibly would be a very useful step toward the
resolution of the question for general graphs.
%
\begin{question}
Is there a deterministic polynomial time $(1+\epsilon)$
approximation algorithm for the supremum of GFF on general graphs?
\end{question}

\begin{remark*}
The question has been solved recently for
general Gaussian process by Meka \cite{Meka12}.
\end{remark*}

\textit{Revisiting the isomorphism theorem.} The isomorphism
theorem was proved by demonstrating an equality of Laplace
transforms for both sides, and very little intuition was provided.
An insightful proof for (\ref{eqlaw}) would be very interesting.
Alternatively, we feel that a proof for the stochastic domination
(\ref{eq-stochastic-domination}) for general graphs (if it is true)
will also shed a good light on understanding the connections between\vadjust{\goodbreak}
local times and GFFs. Plausibly, proving an identity in law is
feasible by showing an equality for Laplace transforms even
without a deep understanding of the processes, while establishing a
stochastic domination seems to require a much deeper insight on the
intrinsic structure of the processes.
%
\begin{question}
Does (\ref{eq-stochastic-domination}) hold for general graphs?
\end{question}
The preceding question, if it is true, will not only shed a good
light on the isomorphism theorem, but also immediately give the
sharp asymptotics as well as an exponential concentration for cover
times.

\textit{A random Eulerian graph model.} We note that
(\ref{eq-random-eulerian}) actually yields a random Eulerian graph
model: given a set of vertices $V$ and nonnegative weights $w_{u,
v}$, take random multiple directed graph such that the multiplicity
of edge $\langle u, v \rangle$ is an independent Poisson variable
with mean $w_{u, v}$, and then re-weighted by a factor of the number
of arborescences contained in the graph (we restricted our space on
Eulerian graphs). More precisely,
\[
\P(\mathcal{G}) \propto\operatorname{ar}(\mathcal{G}) \prod
_{u\neq v} \frac{w_{u, v}^{j_{u, v}}}{j_{u, v}!},
\]
where $j_{u, v}$ denotes the multiplicity of edge $\langle u,
v\rangle$ in $\mathcal{G}$. This random Eulerian model does not seem
to be bizarre in the first place. Also, we believe that an
understanding of this model, in particular on the behavior of the
degrees, could be a useful step for the asymptotics of cover times.

\section*{Acknowledgments}

We thank Yuval Peres, James Lee, Russ Lyons, Jason Miller and Ofer
Zeitouni for helpful discussions, and thank Allan Sly and Tonci
Antunovic for reading the manuscript and valuable suggestions on
exposition, and thank Asaf Nachmias and David Wilson for useful
suggestions on the introduction. Special thanks go to the probability
group at ETH and University Zurich, in particular to David Belius,
Ji\v{r}\'{i} \v{C}ern\'{y}, Alexander Drewitz, Pierre Nolin,
Bal\'{a}zs R\'{a}th, Artem Sapozhnikov and Alain-Sol Sznitman, for a
very careful reading of the paper and for numerous minor corrections.
Finally, we warmly thank two anonymous referees for very careful and
detailed comments.



\printaddresses


\begin{thebibliography}{45}

\bibitem{BR09}
\begin{barticle}[mr]
\bauthor{\bsnm{Addario-Berry},~\bfnm{Louigi}\binits{L.}} \AND
  \bauthor{\bsnm{Reed},~\bfnm{Bruce}\binits{B.}}
(\byear{2009}).
\btitle{Minima in branching random walks}.
\bjournal{Ann. Probab.}
\bvolume{37}
\bpages{1044--1079}.
\bid{doi={10.1214/08-AOP428}, issn={0091-1798}, mr={2537549}}
\bptok{imsref}%
\end{barticle}
\endbibitem

\bibitem{AKS82}
\begin{barticle}[mr]
\bauthor{\bsnm{Ajtai},~\bfnm{M.}\binits{M.}},
  \bauthor{\bsnm{Koml{\'o}s},~\bfnm{J.}\binits{J.}} \AND
  \bauthor{\bsnm{Szemer{\'e}di},~\bfnm{E.}\binits{E.}}
(\byear{1982}).
\btitle{Largest random component of a {$k$}-cube}.
\bjournal{Combinatorica}
\bvolume{2}
\bpages{1--7}.
\bid{doi={10.1007/BF02579276}, issn={0209-9683}, mr={0671140}}
\bptok{imsref}%
\end{barticle}
\endbibitem

\bibitem{AF}
\begin{bmisc}[auto:STB|2013/04/24|11:25:54]
\bauthor{\bsnm{Aldous},~\bfnm{D.}\binits{D.}} \AND
  \bauthor{\bsnm{Fill},~\bfnm{J.}\binits{J.}}
\bhowpublished{Reversible Markov chains and random walks on graphs. Unpublished
  manuscript. Available at \url{http://www.stat.berkeley.edu/\textasciitilde
  aldous/RWG/book.html}}.
\bptok{imsref}%
\end{bmisc}
\endbibitem

\bibitem{Aldous91}
\begin{barticle}[mr]
\bauthor{\bsnm{Aldous},~\bfnm{David~J.}\binits{D.~J.}}
(\byear{1991}).
\btitle{Random walk covering of some special trees}.
\bjournal{J. Math. Anal. Appl.}
\bvolume{157}
\bpages{271--283}.
\bid{doi={10.1016/0022-247X(91)90149-T}, issn={0022-247X}, mr={1109456}}
\bptok{imsref}%
\end{barticle}
\endbibitem

\bibitem{Aldous91b}
\begin{barticle}[mr]
\bauthor{\bsnm{Aldous},~\bfnm{David~J.}\binits{D.~J.}}
(\byear{1991}).
\btitle{Threshold limits for cover times}.
\bjournal{J. Theoret. Probab.}
\bvolume{4}
\bpages{197--211}.
\bid{doi={10.1007/BF01047002}, issn={0894-9840}, mr={1088401}}
\bptok{imsref}%
\end{barticle}
\endbibitem

\bibitem{ABS04}
\begin{barticle}[mr]
\bauthor{\bsnm{Alon},~\bfnm{Noga}\binits{N.}},
  \bauthor{\bsnm{Benjamini},~\bfnm{Itai}\binits{I.}} \AND
  \bauthor{\bsnm{Stacey},~\bfnm{Alan}\binits{A.}}
(\byear{2004}).
\btitle{Percolation on finite graphs and isoperimetric inequalities}.
\bjournal{Ann. Probab.}
\bvolume{32}
\bpages{1727--1745}.
\bid{doi={10.1214/009117904000000414}, issn={0091-1798}, mr={2073175}}
\bptok{imsref}%
\end{barticle}
\endbibitem

\bibitem{BGM10}
\begin{bmisc}[mr]
\bauthor{\bsnm{Benjamini},~\bfnm{Itai}\binits{I.}},
  \bauthor{\bsnm{Gurel-Gurevich},~\bfnm{Ori}\binits{O.}} \AND
  \bauthor{\bsnm{Morris},~\bfnm{Ben}\binits{B.}}
(\byear{2010}).
\bhowpublished{Linear cover time is exponentially unlikely. Preprint. Available
  at \url{http://arxiv.org/abs/1011.3118}}.
\bptok{imsref}%
\end{bmisc}
\endbibitem

\bibitem{BNP09}
\begin{barticle}[mr]
\bauthor{\bsnm{Benjamini},~\bfnm{Itai}\binits{I.}},
  \bauthor{\bsnm{Nachmias},~\bfnm{Asaf}\binits{A.}} \AND
  \bauthor{\bsnm{Peres},~\bfnm{Yuval}\binits{Y.}}
(\byear{2011}).
\btitle{Is the critical percolation probability local?}
\bjournal{Probab. Theory Related Fields}
\bvolume{149}
\bpages{261--269}.
\bid{doi={10.1007/s00440-009-0251-5}, issn={0178-8051}, mr={2773031}}
\bptnote{check year}%
\bptok{imsref}%
\end{barticle}
\endbibitem

\bibitem{BDG01}
\begin{barticle}[mr]
\bauthor{\bsnm{Bolthausen},~\bfnm{Erwin}\binits{E.}},
  \bauthor{\bsnm{Deuschel},~\bfnm{Jean-Dominique}\binits{J.-D.}} \AND
  \bauthor{\bsnm{Giacomin},~\bfnm{Giambattista}\binits{G.}}
(\byear{2001}).
\btitle{Entropic repulsion and the maximum of the two-dimensional harmonic
  crystal}.
\bjournal{Ann. Probab.}
\bvolume{29}
\bpages{1670--1692}.
\bid{doi={10.1214/aop/1015345767}, issn={0091-1798}, mr={1880237}}
\bptok{imsref}%
\end{barticle}
\endbibitem

\bibitem{BDZ10}
\begin{barticle}[mr]
\bauthor{\bsnm{Bolthausen},~\bfnm{Erwin}\binits{E.}},
  \bauthor{\bsnm{Deuschel},~\bfnm{Jean~Dominique}\binits{J.~D.}} \AND
  \bauthor{\bsnm{Zeitouni},~\bfnm{Ofer}\binits{O.}}
(\byear{2011}).
\btitle{Recursions and tightness for the maximum of the discrete, two
  dimensional {G}aussian free field}.
\bjournal{Electron. Commun. Probab.}
\bvolume{16}
\bpages{114--119}.
\bid{doi={10.1214/ECP.v16-1610}, issn={1083-589X}, mr={2772390}}
\bptok{imsref}%
\end{barticle}
\endbibitem

\bibitem{BZ09}
\begin{barticle}[mr]
\bauthor{\bsnm{Bramson},~\bfnm{Maury}\binits{M.}} \AND
  \bauthor{\bsnm{Zeitouni},~\bfnm{Ofer}\binits{O.}}
(\byear{2009}).
\btitle{Tightness for a family of recursion equations}.
\bjournal{Ann. Probab.}
\bvolume{37}
\bpages{615--653}.
\bid{doi={10.1214/08-AOP414}, issn={0091-1798}, mr={2510018}}
\bptok{imsref}%
\end{barticle}
\endbibitem

\bibitem{BZ10}
\begin{barticle}[mr]
\bauthor{\bsnm{Bramson},~\bfnm{Maury}\binits{M.}} \AND
  \bauthor{\bsnm{Zeitouni},~\bfnm{Ofer}\binits{O.}}
(\byear{2012}).
\btitle{Tightness of the recentered maximum of the two-dimensional discrete
  {G}aussian free field}.
\bjournal{Comm. Pure Appl. Math.}
\bvolume{65}
\bpages{1--20}.
\bid{doi={10.1002/cpa.20390}, issn={0010-3640}, mr={2846636}}
\bptok{imsref}%
\end{barticle}
\endbibitem

\bibitem{Bramson78}
\begin{barticle}[mr]
\bauthor{\bsnm{Bramson},~\bfnm{Maury~D.}\binits{M.~D.}}
(\byear{1978}).
\btitle{Maximal displacement of branching {B}rownian motion}.
\bjournal{Comm. Pure Appl. Math.}
\bvolume{31}
\bpages{531--581}.
\bid{issn={0010-3640}, mr={0494541}}
\bptok{imsref}%
\end{barticle}
\endbibitem

\bibitem{CRRST96}
\begin{barticle}[mr]
\bauthor{\bsnm{Chandra},~\bfnm{Ashok~K.}\binits{A.~K.}},
  \bauthor{\bsnm{Raghavan},~\bfnm{Prabhakar}\binits{P.}},
  \bauthor{\bsnm{Ruzzo},~\bfnm{Walter~L.}\binits{W.~L.}},
  \bauthor{\bsnm{Smolensky},~\bfnm{Roman}\binits{R.}} \AND
  \bauthor{\bsnm{Tiwari},~\bfnm{Prasoon}\binits{P.}}
(\byear{1996/97}).
\btitle{The electrical resistance of a graph captures its commute and cover
  times}.
\bjournal{Comput. Complexity}
\bvolume{6}
\bpages{312--340}.
\bid{doi={10.1007/BF01270385}, issn={1016-3328}, mr={1613611}}
\bptok{imsref}%
\end{barticle}
\endbibitem

\bibitem{Chatterjee08}
\begin{bmisc}[auto:STB|2013/04/24|11:25:54]
\bauthor{\bsnm{Chatterjee},~\bfnm{S.}\binits{S.}}
(\byear{2008}).
\bhowpublished{Chaos, concentration, and multiple valleys. Preprint. Available
  at \url{http://arxiv.org/abs/0810.4221}}.
\bptok{imsref}%
\end{bmisc}
\endbibitem

\bibitem{DPRZ04}
\begin{barticle}[mr]
\bauthor{\bsnm{Dembo},~\bfnm{Amir}\binits{A.}},
  \bauthor{\bsnm{Peres},~\bfnm{Yuval}\binits{Y.}},
  \bauthor{\bsnm{Rosen},~\bfnm{Jay}\binits{J.}} \AND
  \bauthor{\bsnm{Zeitouni},~\bfnm{Ofer}\binits{O.}}
(\byear{2004}).
\btitle{Cover times for {B}rownian motion and random walks in two dimensions}.
\bjournal{Ann. of Math. (2)}
\bvolume{160}
\bpages{433--464}.
\bid{doi={10.4007/annals.2004.160.433}, issn={0003-486X}, mr={2123929}}
\bptok{imsref}%
\end{barticle}
\endbibitem

\bibitem{DLP10}
\begin{barticle}[mr]
\bauthor{\bsnm{Ding},~\bfnm{Jian}\binits{J.}},
  \bauthor{\bsnm{Lee},~\bfnm{James~R.}\binits{J.~R.}} \AND
  \bauthor{\bsnm{Peres},~\bfnm{Yuval}\binits{Y.}}
(\byear{2012}).
\btitle{Cover times, blanket times, and majorizing measures}.
\bjournal{Ann. of Math. (2)}
\bvolume{175}
\bpages{1409--1471}.
\bid{doi={10.4007/annals.2012.175.3.8}, issn={0003-486X}, mr={2912708}}
\bptok{imsref}%
\end{barticle}
\endbibitem

\bibitem{DZ11}
\begin{barticle}[mr]
\bauthor{\bsnm{Ding},~\bfnm{Jian}\binits{J.}} \AND
  \bauthor{\bsnm{Zeitouni},~\bfnm{Ofer}\binits{O.}}
(\byear{2012}).
\btitle{A sharp estimate for cover times on binary trees}.
\bjournal{Stochastic Process. Appl.}
\bvolume{122}
\bpages{2117--2133}.
\bid{doi={10.1016/j.spa.2012.03.008}, issn={0304-4149}, mr={2921974}}
\bptok{imsref}%
\end{barticle}
\endbibitem

\bibitem{Durrett}
\begin{bbook}[mr]
\bauthor{\bsnm{Durrett},~\bfnm{Rick}\binits{R.}}
(\byear{2010}).
\btitle{Probability: Theory and Examples},
\bedition{4th} ed.
\bpublisher{Cambridge Univ. Press}, \blocation{Cambridge}.
\bid{mr={2722836}}
\bptok{imsref}%
\end{bbook}
\endbibitem

\bibitem{Dynkin80}
\begin{barticle}[mr]
\bauthor{\bsnm{Dynkin},~\bfnm{E.~B.}\binits{E.~B.}}
(\byear{1980}).
\btitle{Markov processes and random fields}.
\bjournal{Bull. Amer. Math. Soc. (N.S.)}
\bvolume{3}
\bpages{975--999}.
\bid{doi={10.1090/S0273-0979-1980-14831-4}, issn={0273-0979}, mr={0585179}}
\bptok{imsref}%
\end{barticle}
\endbibitem

\bibitem{Dynkin84}
\begin{barticle}[mr]
\bauthor{\bsnm{Dynkin},~\bfnm{E.~B.}\binits{E.~B.}}
(\byear{1984}).
\btitle{Gaussian and non-{G}aussian random fields associated with {M}arkov
  processes}.
\bjournal{J. Funct. Anal.}
\bvolume{55}
\bpages{344--376}.
\bid{doi={10.1016/0022-1236(84)90004-1}, issn={0022-1236}, mr={0734803}}
\bptok{imsref}%
\end{barticle}
\endbibitem

\bibitem{Dynkin83}
\begin{bincollection}[mr]
\bauthor{\bsnm{Dynkin},~\bfnm{E.~B.}\binits{E.~B.}}
(\byear{1984}).
\btitle{Local times and quantum fields}.
In \bbooktitle{Seminar on Stochastic Processes, 1983 ({G}ainesville, {F}la.,
  1983)}.
\bseries{Progr. Probab. Statist.}
\bvolume{7}
\bpages{69--83}.
\bpublisher{Birkh\"auser}, \blocation{Boston, MA}.
\bid{mr={0902412}}
\bptok{imsref}%
\end{bincollection}
\endbibitem

\bibitem{Eisenbaum95}
\begin{bincollection}[mr]
\bauthor{\bsnm{Eisenbaum},~\bfnm{Nathalie}\binits{N.}}
(\byear{1995}).
\btitle{Une version sans conditionnement du th\'eor\`eme d'isomorphisms de
  {D}ynkin}.
In \bbooktitle{S\'eminaire de {P}robabilit\'es, {XXIX}}.
\bseries{Lecture Notes in Math.}
\bvolume{1613}
\bpages{266--289}.
\bpublisher{Springer}, \blocation{Berlin}.
\bid{doi={10.1007/BFb0094219}, mr={1459468}}
\bptok{imsref}%
\end{bincollection}
\endbibitem

\bibitem{EKMRS00}
\begin{barticle}[mr]
\bauthor{\bsnm{Eisenbaum},~\bfnm{Nathalie}\binits{N.}},
  \bauthor{\bsnm{Kaspi},~\bfnm{Haya}\binits{H.}},
  \bauthor{\bsnm{Marcus},~\bfnm{Michael~B.}\binits{M.~B.}},
  \bauthor{\bsnm{Rosen},~\bfnm{Jay}\binits{J.}} \AND
  \bauthor{\bsnm{Shi},~\bfnm{Zhan}\binits{Z.}}
(\byear{2000}).
\btitle{A {R}ay--{K}night theorem for symmetric {M}arkov processes}.
\bjournal{Ann. Probab.}
\bvolume{28}
\bpages{1781--1796}.
\bid{doi={10.1214/aop/1019160507}, issn={0091-1798}, mr={1813843}}
\bptok{imsref}%
\end{barticle}
\endbibitem

\bibitem{FZ09}
\begin{bmisc}[auto:STB|2013/04/24|11:25:54]
\bauthor{\bsnm{Feige},~\bfnm{U.}\binits{U.}} \AND
  \bauthor{\bsnm{Zeitouni},~\bfnm{O.}\binits{O.}}
(\byear{2009}).
\bhowpublished{Deterministic approximation for the cover time of trees.
  Preprint. Available at \url{http://arxiv.org/abs/0909.2005v1}}.
\bptok{imsref}%
\end{bmisc}
\endbibitem

\bibitem{Fernique74}
\begin{barticle}[mr]
\bauthor{\bsnm{Fernique},~\bfnm{Xavier}\binits{X.}}
(\byear{1974}).
\btitle{Des r\'esultats nouveaux sur les processus gaussiens}.
\bjournal{C. R. Acad. Sci. Paris S\'er. A}
\bvolume{278}
\bpages{363--365}.
\bid{mr={0413237}}
\bptok{imsref}%
\end{barticle}
\endbibitem

\bibitem{Janson97}
\begin{bbook}[mr]
\bauthor{\bsnm{Janson},~\bfnm{Svante}\binits{S.}}
(\byear{1997}).
\btitle{Gaussian {H}ilbert Spaces}.
\bseries{Cambridge Tracts in Mathematics}
\bvolume{129}.
\bpublisher{Cambridge Univ. Press}, \blocation{Cambridge}.
\bid{doi={10.1017/CBO9780511526169}, mr={1474726}}
\bptok{imsref}%
\end{bbook}
\endbibitem

\bibitem{KKLV00}
\begin{bincollection}[mr]
\bauthor{\bsnm{Kahn},~\bfnm{J.}\binits{J.}},
  \bauthor{\bsnm{Kim},~\bfnm{J.~H.}\binits{J.~H.}},
  \bauthor{\bsnm{Lov{\'a}sz},~\bfnm{L.}\binits{L.}} \AND
  \bauthor{\bsnm{Vu},~\bfnm{V.~H.}\binits{V.~H.}}
(\byear{2000}).
\btitle{The cover time, the blanket time, and the {M}atthews bound}.
In \bbooktitle{41st {A}nnual {S}ymposium on {F}oundations of {C}omputer
  {S}cience ({R}edondo {B}each, {CA}, 2000)}
\bpages{467--475}.
\bpublisher{IEEE Comput. Soc.}, \blocation{Los Alamitos, CA}.
\bid{doi={10.1109/SFCS.2000.892134}, mr={1931843}}
\bptok{imsref}%
\end{bincollection}
\endbibitem

\bibitem{KR93}
\begin{barticle}[mr]
\bauthor{\bsnm{Klein},~\bfnm{D.~J.}\binits{D.~J.}} \AND
  \bauthor{\bsnm{Randi{\'c}},~\bfnm{M.}\binits{M.}}
(\byear{1993}).
\btitle{Resistance distance}.
\bjournal{J. Math. Chem.}
\bvolume{12}
\bpages{81--95}.
\bnote{Applied graph theory and discrete mathematics in chemistry (Saskatoon,
  SK, 1991)}.
\bid{doi={10.1007/BF01164627}, issn={0259-9791}, mr={1219566}}
\bptok{imsref}%
\end{barticle}
\endbibitem

\bibitem{Knight63}
\begin{barticle}[mr]
\bauthor{\bsnm{Knight},~\bfnm{F.~B.}\binits{F.~B.}}
(\byear{1963}).
\btitle{Random walks and a sojourn density process of {B}rownian motion}.
\bjournal{Trans. Amer. Math. Soc.}
\bvolume{109}
\bpages{56--86}.
\bid{issn={0002-9947}, mr={0154337}}
\bptok{imsref}%
\end{barticle}
\endbibitem

\bibitem{Ledoux89}
\begin{bbook}[mr]
\bauthor{\bsnm{Ledoux},~\bfnm{Michel}\binits{M.}}
(\byear{2001}).
\btitle{The Concentration of Measure Phenomenon}.
\bseries{Mathematical Surveys and Monographs}
\bvolume{89}.
\bpublisher{Amer. Math. Soc.}, \blocation{Providence, RI}.
\bid{mr={1849347}}
\bptok{imsref}%
\end{bbook}
\endbibitem

\bibitem{LPW09}
\begin{bbook}[mr]
\bauthor{\bsnm{Levin},~\bfnm{David~A.}\binits{D.~A.}},
  \bauthor{\bsnm{Peres},~\bfnm{Yuval}\binits{Y.}} \AND
  \bauthor{\bsnm{Wilmer},~\bfnm{Elizabeth~L.}\binits{E.~L.}}
(\byear{2009}).
\btitle{Markov Chains and Mixing Times}.
\bpublisher{Amer. Math. Soc.}, \blocation{Providence, RI}.
\bnote{With a chapter by James G. Propp and David B. Wilson}.
\bid{mr={2466937}}
\bptok{imsref}%
\end{bbook}
\endbibitem

\bibitem{Lov96}
\begin{bincollection}[mr]
\bauthor{\bsnm{Lov{\'a}sz},~\bfnm{L.}\binits{L.}}
(\byear{1996}).
\btitle{Random walks on graphs: A survey}.
In \bbooktitle{Combinatorics, {P}aul {E}rd{\H o}s Is Eighty, {V}ol.\ 2
  ({K}eszthely, 1993)}.
\bseries{Bolyai Society Mathematical Studies}
\bvolume{2}
\bpages{353--397}.
\bpublisher{J\'anos Bolyai Math. Soc.}, \blocation{Budapest}.
\bid{mr={1395866}}
\bptok{imsref}%
\end{bincollection}
\endbibitem

\bibitem{LP}
\begin{bmisc}[auto:STB|2013/04/24|11:25:54]
\bauthor{\bsnm{Lyons},~\bfnm{R.}\binits{R.}} \AND
  \bauthor{\bsnm{Peres},~\bfnm{Y.}\binits{Y.}}
(\byear{2009}).
\bhowpublished{Probability on trees and networks. Unpublished manuscript.
  Current version available at \url{http://mypage.iu.edu/\textasciitilde
  rdlyons/prbtree/book.pdf}}.
\bptok{imsref}%
\end{bmisc}
\endbibitem

\bibitem{MR92}
\begin{barticle}[mr]
\bauthor{\bsnm{Marcus},~\bfnm{Michael~B.}\binits{M.~B.}} \AND
  \bauthor{\bsnm{Rosen},~\bfnm{Jay}\binits{J.}}
(\byear{1992}).
\btitle{Sample path properties of the local times of strongly symmetric
  {M}arkov processes via {G}aussian processes}.
\bjournal{Ann. Probab.}
\bvolume{20}
\bpages{1603--1684}.
\bid{issn={0091-1798}, mr={1188037}}
\bptok{imsref}%
\end{barticle}
\endbibitem

\bibitem{MR01}
\begin{bincollection}[mr]
\bauthor{\bsnm{Marcus},~\bfnm{Michael~B.}\binits{M.~B.}} \AND
  \bauthor{\bsnm{Rosen},~\bfnm{Jay}\binits{J.}}
(\byear{2001}).
\btitle{Gaussian processes and local times of symmetric {L}\'evy processes}.
In \bbooktitle{L\'evy Processes}
\bpages{67--88}.
\bpublisher{Birkh\"auser}, \blocation{Boston, MA}.
\bid{mr={1833693}}
\bptok{imsref}%
\end{bincollection}
\endbibitem

\bibitem{MR06}
\begin{bbook}[mr]
\bauthor{\bsnm{Marcus},~\bfnm{Michael~B.}\binits{M.~B.}} \AND
  \bauthor{\bsnm{Rosen},~\bfnm{Jay}\binits{J.}}
(\byear{2006}).
\btitle{Markov Processes, {G}aussian Processes, and Local Times}.
\bseries{Cambridge Studies in Advanced Mathematics}
\bvolume{100}.
\bpublisher{Cambridge Univ. Press}, \blocation{Cambridge}.
\bid{doi={10.1017/CBO9780511617997}, mr={2250510}}
\bptok{imsref}%
\end{bbook}
\endbibitem

\bibitem{Meka12}
\begin{bmisc}[auto:STB|2013/04/24|11:25:54]
\bauthor{\bsnm{Meka},~\bfnm{R.}\binits{R.}}
(\byear{2012}).
\bhowpublished{A PTAS for computing the supremum of Gaussian processes. In
\textit{Proceedings of the 2012 IEEE 53rd Annual Symposium on Foundations of Computer Science} 217--222.
IEEE Computer Society, Washington, DC}.
\bptok{imsref}%
\end{bmisc}
\endbibitem

\bibitem{MP09}
\begin{barticle}[mr]
\bauthor{\bsnm{Miller},~\bfnm{Jason}\binits{J.}} \AND
  \bauthor{\bsnm{Peres},~\bfnm{Yuval}\binits{Y.}}
(\byear{2012}).
\btitle{Uniformity of the uncovered set of random walk and cutoff for
  lamplighter chains}.
\bjournal{Ann. Probab.}
\bvolume{40}
\bpages{535--577}.
\bid{doi={10.1214/10-AOP624}, issn={0091-1798}, mr={2952084}}
\bptok{imsref}%
\end{barticle}
\endbibitem

\bibitem{Ray63}
\begin{barticle}[mr]
\bauthor{\bsnm{Ray},~\bfnm{Daniel}\binits{D.}}
(\byear{1963}).
\btitle{Sojourn times of diffusion processes}.
\bjournal{Illinois J. Math.}
\bvolume{7}
\bpages{615--630}.
\bid{issn={0019-2082}, mr={0156383}}
\bptok{imsref}%
\end{barticle}
\endbibitem

\bibitem{Talagrand87}
\begin{barticle}[mr]
\bauthor{\bsnm{Talagrand},~\bfnm{Michel}\binits{M.}}
(\byear{1987}).
\btitle{Regularity of {G}aussian processes}.
\bjournal{Acta Math.}
\bvolume{159}
\bpages{99--149}.
\bid{doi={10.1007/BF02392556}, issn={0001-5962}, mr={0906527}}
\bptok{imsref}%
\end{barticle}
\endbibitem

\bibitem{Talagrand05}
\begin{bbook}[mr]
\bauthor{\bsnm{Talagrand},~\bfnm{Michel}\binits{M.}}
(\byear{2005}).
\btitle{The Generic Chaining: Upper and Lower Bounds of Stochastic Processes}.
\bpublisher{Springer}, \blocation{Berlin}.
\bid{mr={2133757}}
\bptok{imsref}%
\end{bbook}
\endbibitem

\bibitem{TS41}
\begin{barticle}[mr]
\bauthor{\bsnm{Tutte},~\bfnm{W.~T.}\binits{W.~T.}} \AND
  \bauthor{\bsnm{Smith},~\bfnm{C.~A.~B.}\binits{C.~A.~B.}}
(\byear{1941}).
\btitle{On {u}nicursal {p}aths in a {n}etwork of {d}egree 4}.
\bjournal{Amer. Math. Monthly}
\bvolume{48}
\bpages{233--237}.
\bid{doi={10.2307/2302716}, issn={0002-9890}, mr={1525117}}
\bptok{imsref}%
\end{barticle}
\endbibitem

\bibitem{BE51}
\begin{barticle}[mr]
\bauthor{\bparticle{van} \bsnm{Aardenne-Ehrenfest},~\bfnm{T.}\binits{T.}} \AND
  \bauthor{\bparticle{de} \bsnm{Bruijn},~\bfnm{N.~G.}\binits{N.~G.}}
(\byear{1951}).
\btitle{Circuits and trees in oriented linear graphs}.
\bjournal{Simon Stevin}
\bvolume{28}
\bpages{203--217}.
\bid{issn={0037-5454}, mr={0047311}}
\bptok{imsref}%
\end{barticle}
\endbibitem

\end{thebibliography}
\end{document}